\documentclass[a4paper,12pt]{article}
\usepackage{amsmath}
\usepackage{amsfonts}
\usepackage{amssymb}
\usepackage{latexsym}
\usepackage{epsfig}
\usepackage{graphicx}
\usepackage{oldgerm}
\usepackage{theorem}

\def\C{\mathbb{C}}
\def\R{\mathbb{R}}
\def\N{\mathbb{N}}
\def\W{\mathbb{W}}
\def\Z{\mathbb{Z}}

\def\P{{\mathbb P}}
\def\E{{\mathbb E}}

\def\mbK{\mathbb{K}}

\def\bK{{\bf K}}
\def\bP{{\bf P}}

\def\x{\mib{x}}
\def\y{\mib{y}}
\def\bB{\mib{B}}

\def\vchi{\mib{\chi}}

\def\mM{\mathfrak{M}}
\def\mX{\mathfrak{X}}

\def\supp{{\rm supp}\ }
\def\cG{{\cal G}}

\def\Det{\mathop{\mathrm{Det}}}

\theorembodyfont{\itshape}

\newtheorem{thm}{Theorem}[section]
\newtheorem{lem}[thm]{Lemma}

\newcommand{\mib}[1]{\mbox{\boldmath $#1$}}
\newcommand{\SSC}[1]{\section{#1}\setcounter{equation}{0}}
\newcommand{\qed}{\hbox{\rule[-2pt]{3pt}{6pt}}}

\begin{document}

\title{\bf 
Noncolliding Squared Bessel Processes}
\author{
Makoto Katori
\footnote{
Department of Physics,
Faculty of Science and Engineering,
Chuo University, 
Kasuga, Bunkyo-ku, Tokyo 112-8551, Japan;
e-mail: katori@phys.chuo-u.ac.jp
}
and 
Hideki Tanemura
\footnote{
Department of Mathematics and Informatics,
Faculty of Science, Chiba University, 
1-33 Yayoi-cho, Inage-ku, Chiba 263-8522, Japan;
e-mail: tanemura@math.s.chiba-u.ac.jp
}}
\date{16 January 2011}
\pagestyle{plain}
\maketitle
\begin{abstract}
We consider a particle system of the squared Bessel processes
with index $\nu > -1$ conditioned never to collide with
each other,
in which if $-1 < \nu < 0$ the origin is 
assumed to be reflecting.
When the number of particles is finite, we prove 
for any fixed initial configuration that 
this noncolliding diffusion process is determinantal
in the sense that any multitime correlation function is
given by a determinant with a continuous kernel called
the correlation kernel.
When the number of particles is infinite,
we give sufficient conditions for initial configurations
so that the system is well defined.
There the process with an infinite number of particles
is determinantal and the correlation kernel is
expressed using an entire function
represented by the Weierstrass canonical product,
whose zeros on the positive part of the real axis 
are given by the particle-positions in the initial
configuration.
From the class of infinite-particle initial configurations
satisfying our conditions, we report one example in detail, 
which is a fixed configuration such that
every point of the square of positive zero of the
Bessel function $J_{\nu}$ is occupied by one particle. 
The process starting from this initial
configuration shows a relaxation phenomenon
converging to the stationary process,
which is determinantal with the extended
Bessel kernel, in the long-term limit.
\end{abstract}

\normalsize

\SSC{Introduction}

Let $\mM$ be the space of nonnegative integer-valued
Radon measures on $\R$,
which is a Polish space with the vague topology.
We say $\xi_n \in \mM, n \in \N \equiv \{1,2, \dots\}$ 
converges to $\xi \in \mM$ weakly in the vague topology, if 
$\lim_{n \to \infty} \int_{\R} \varphi(x) \xi_n(dx)
=\int_{\R} \varphi(x) \xi(dx)$ 
for any $\varphi \in {\rm C}_0(\R)$,
where ${\rm C}_0(\R)$ is the set of all 
continuous real-valued functions with
compact supports in $\R$.
Any element of $\mM$ can be represented by 
$\sum_{i \in \Lambda} \delta_{x_i}(\cdot)$
with a sequence of points in $\R$,
$\x=(x_i)_{i \in \Lambda}$, satisfying
$\sharp\{x_i: x_i \in I \} < \infty$ for any
compact subset $I \subset \R$,
and with a countable set (an index set) $\Lambda$.
We call an element $\xi$ of $\mM$ an unlabeled configuration,
and a sequence of points $\x$ a labeled configuration.
For $A \subset \R$, we write 
the restriction of $\xi \in \mM$ on $A$ as 
$(\xi\cap A) (\cdot)=\sum_{i \in \Lambda : x_i \in  A}
\delta_{x_i}(\cdot)$.
Let $\R_+ =\{x \in \R: x \geq 0 \}$
and define 
$\mM^+=\{(\xi \cap \R_+)(\cdot) : \xi(\cdot) \in \mM\}$.
In the present paper
we consider a one-parameter family
of $\mM^+$-valued processes 
with a parameter $\nu > -1$,
\begin{equation}
\Xi^{(\nu)}(t, \cdot)=\sum_i \delta_{X^{(\nu)}_i(t)}(\cdot),
\quad t \in [0, \infty),
\label{eqn:Xi1}
\end{equation}
describing a particle system of squared Bessel processes
with index $\nu > -1$ (BESQ$^{(\nu)}$) 
interacting with each other by {\it long-ranged repulsive forces},
such that $X^{(\nu)}_i(t)$'s satisfy the SDEs
\begin{eqnarray}
d X^{(\nu)}_i(t) &=& 2 \sqrt{X^{(\nu)}_i(t)} dB_i(t)
+2 (\nu+1) dt \nonumber\\
&+& 4 X^{(\nu)}_i(t) \sum_{j: j \not=i}
\frac{1}{X^{(\nu)}_i(t)-X^{(\nu)}_j(t)} dt, \quad
i=1,2, \cdots, \quad t \in [0, \infty)
\label{eqn:ncBESQ}
\end{eqnarray}
with a collection of independent standard
Brownian motions (BMs), 
$\{B_i(t), i \in \N\}$,
and, if $-1 < \nu < 0$, with a reflection wall at the origin.
Note that for the BM in $\R^{d}$,
$\widetilde{\bB}(t)=(\widetilde{B}_1(t), \dots, \widetilde{B}_d(t)), 
d \in \N$,
the square of its distance from the origin,
$X(t) \equiv |\widetilde{\bB}(t)|^2
=\sum_{i=1}^{d} \widetilde{B}_i(t)^2$, solves the SDE,
$dX(t)=2 \sqrt{X(t)} dB(t) + 2(\nu+1) dt$
with
\begin{equation}
  \nu=\frac{d}{2}-1,
\label{eqn:dnu}
\end{equation}
where $B(t)$ is a standard BM which is different from
$\widetilde{B}_{i}(t), 1 \leq i \leq d$ \cite{RY98,BS02}.
We give the initial configuration of the process
$\xi(\cdot)=\Xi^{(\nu)}(0, \cdot)
=\sum_{i} \delta_{x_i}(\cdot)$
and the process is denoted by $(\Xi^{(\nu)}(t), \P^{\xi}_{\nu})$.

When the number of particles is finite,
$\xi_N(\R_+)=N < \infty$,
the process $(\Xi^{(\nu)}(t), \P^{\xi_N}_{\nu})$ is realized
in the following systems.
\begin{description}
\item{(i)} \quad
When $\nu \in \N_0 \equiv \N \cup \{0\}$,
{\it i.e.}, when 
the corresponding dimension $d$ given by (\ref{eqn:dnu})
is a positive even integer, 
$(\Xi^{(\nu)}(t), \P^{\xi_N}_{\nu})$ is realized as the
eigenvalue process of the
{\it Laguerre process} \cite{KO01}.
Let $M(t)$ be an $(N+\nu) \times N$
matrix, whose entries are independent 
complex BMs having the real and imaginary
parts given by independent standard BMs, and set
$L(t)=M(t)^* M(t)$.
The $N \times N$ matrix-valued process
${\rm L}=(L(t))_{t \in [0, \infty)}$
is called the Laguerre process.
The matrix $L(t)$ is Hermitian and positive definite,
and its $N$ eigenvalues satisfy (\ref{eqn:ncBESQ})
with $i=1,2, \dots,N$.
When the entries of $M(t)$ are independent standard
{\it real} BMs, the matrix-valued process
$(M(t)^{\rm T} M(t))_{t \in [0, \infty)}$ is called
the Wishart process \cite{Bru91}, and thus
${\rm L}$ is also called the
{\it complex Wishart process}.
The eigenvalue processes of the real and complex
Wishart processes are related with the random matrix
theory \cite{Meh04,For10} 
for the {\it chiral Gaussian ensembles}
studied in the high energy physics 
(see \cite{KT04} and references therein).

\item{(ii)} \quad
Let ${\cal H}(N)$ be the space of $N \times N$
Hermitian matrices. And let
$\mathfrak{sp}(2N, \C)$ and $\mathfrak{so}(2N, \C)$
be the symplectic Lie algebra and
the orthogonal Lie algebra,
having $2N \times 2N$-matrix representations,
respectively.
If the $2N \times 2N$ matrix is in the space
${\cal H}_{\rm C}(2N) \equiv {\cal H}(2N) \cap
\mathfrak{sp}(2N, \C)$ or in 
${\cal H}_{\rm D}(2N) \equiv {\cal H}(2N) \cap
\mathfrak{so}(2N, \C)$,
its eigenvalues are given by $N$ pairs of 
positive and negative ones with the same absolute value,
$\{(\lambda_i, - \lambda_i): \lambda_i \geq 0,
1 \leq i \leq N\}$.
Consider the ${\cal H}_{\rm C}(2N)$-valued
and the ${\cal H}_{\rm D}(2N)$-valued
Brownian motions.
The dynamics of positive eigenvalues
of them are described by (\ref{eqn:ncBESQ})
with $\nu=1/2 \, (d=3)$ for the former case
and with $\nu=-1/2 \, (d=1)$ for the latter case,
respectively \cite{KT04}.
The pairing of positive and negative eigenvalues
simulates the particle-hole symmetry 
in the energy space of the Bogoliubov-de Gennes formalism of
superconductivity and these processes 
are related with the random matrix theory
studied in the solid-state physics \cite{AZ97}.

\item{(iii)} \quad
Let $p^{(\nu)}(t, y|x), y \in \R_+$,
be the transition probability density for BESQ$^{(\nu)}$,
$\nu > -1$,
\begin{equation}
p^{(\nu)}(t, y|x)
=  \left\{ \begin{array}{ll}
\displaystyle{
\frac{1}{2t} \left( \frac{y}{x} \right)^{\nu/2}
\exp \left( - \frac{x+y}{2t} \right)
I_{\nu} \left( \frac{\sqrt{xy}}{t} \right)},
& \quad t>0, x>0, \cr
\displaystyle{
\frac{y^{\nu}}{(2t)^{\nu+1} \Gamma(\nu+1)} e^{-y/2t}},
& \quad t >0, x=0, \cr
& \cr
\delta(y-x),
& \quad t=0, x \in \R_+,
\end{array} \right.
\label{eqn:p0}
\end{equation}
if $-1 < \nu < 0$, the origin is assumed to be reflecting
\cite{RY98,BS02},
where $I_{\nu}(x)$ is the modified Bessel function
of the first kind defined by
\begin{equation}
I_{\nu}(x) = \sum_{n=0}^{\infty} 
\frac{1}{\Gamma(n+1) \Gamma(n+1+\nu)}
\left( \frac{x}{2} \right)^{2n+\nu}
\label{eqn:I1}
\end{equation}
with the Gamma function
$\Gamma(z) = \int_{0}^{\infty} e^{-u} u^{z-1} du, 
\ \Re u > 0$.
First we consider the following Karlin-McGregor determinant 
\cite{KM59},
\begin{equation}
f_N^{(\nu)}(t, \y|\x)
= \det_{1 \leq i, j \leq N}
\Big[ p^{(\nu)}(t, y_i|x_j) \Big],
\quad t \geq 0, \quad \x=(x_i)_{i=1}^{N}, 
\y=(y_i)_{i=1}^{N} \in \W_N^{+},
\label{eqn:fN}
\end{equation}
where $\W_N^{+}=\{\x=(x_1, \dots, x_N) \in \R^N:
0 \leq x_1<\cdots <x_N\}$.
The transition probability density of the $N$-particle
system of BESQ$^{(\nu)}$ conditioned never to collide
with each other, which we call the {\it noncolliding BESQ}$^{(\nu)}$,
is given by the
$h$-transform of (\ref{eqn:fN}),
\begin{equation}
p_N^{(\nu)}(t, \y|\x)
=h_N(\y) f_N^{(\nu)}(t, \y|\x) \frac{1}{h_N(\x)}
\label{eqn:pNnu}
\end{equation}
with the harmonic function given by
the Vandermonde determinant \cite{Gra99,KO01,KT04}
\begin{equation}
h_N(\x)=\det_{1 \leq i, j \leq N} [x_i^{j-1}]
=\prod_{1 \leq i < j \leq N} (x_j-x_i).
\label{eqn:h1}
\end{equation}
It is easy to confirm that $p_N^{(\nu)}(t, \cdot|\x)$ 
satisfies the following backward Kolmogorov equation
\begin{eqnarray}
\frac{\partial}{\partial t} u(t,\x)
&=& 2 \sum_{i=1}^{N} x_i 
\frac{\partial^2}{\partial x_i^2} u(t,\x)
+2 (\nu+1) \sum_{i=1}^{N} \frac{\partial}{\partial x_i} u(t,\x)
\nonumber\\
&+& 4 \sum_{i=1}^{N} \sum_{\substack{j=1 \\ j \not= i}}^{N}
\frac{x_i}{x_i-x_j} \frac{\partial}{\partial x_i} u(t,\x),
\quad t \geq 0, \quad \x \in \W_N^{+},
\nonumber
\end{eqnarray}
and it implies that the process 
$(\Xi^{(\nu)}(t), \P^{\xi_N}_{\nu}), \nu > -1$,
is realized as the noncolliding BESQ$^{(\nu)}$ \cite{KT04}.
We put
\begin{equation}
\mM^+_0= \{ \xi\in\mM^+ : \xi(\{x\})\le 1 \mbox { for any }  x\in\R
\}.
\label{eqn:mM0}
\end{equation}
We see that $\Xi^{(\nu)}(t) \in \mM^+_0, \forall t > 0$.
\end{description}

In the present paper, we call the process
$(\Xi^{(\nu)}(t), \P^{\xi}_{\nu})$ 
the noncolliding BESQ$^{(\nu)}$.
See \cite{TW07,KT07a,KMW09} for related
noncolliding diffusion processes.

Assume that $\xi_N \in \mM^+$ with $\xi_N(\R_+)=N \in \N$.
For any $M \in \N$ and any time sequence
$0 < t_1 < \cdots < t_M < \infty$, 
the formula (\ref{eqn:pNnu}) and the Markov property of
the system give the {\it multitime probability density}
of $(\Xi^{(\nu)}, \P^{\xi_N}_{\nu})$ as \cite{Gra99,KT07b}
\begin{equation}
p_{\nu}^{\xi_N}(t_1, \xi_N^{(1)}; \dots ; t_M, \xi_N^{(M)})
= h_N(\x^{(M)}) \prod_{m=1}^{M-1} 
f_N^{(\nu)}(t_{m+1}-t_{m}; \x^{(m+1)}|\x^{(m)})
\frac{f_N^{(\nu)}(t_1, \x^{(1)}|\x)}{h_N(\x)}
\label{eqn:pxi}
\end{equation}
with $\xi_N(\cdot)=\sum_{i=1}^{N} \delta_{x_i}(\cdot) \in \mM^+$,
$0 \leq x_1 \leq x_2 \leq \cdots
\leq x_N$ for the initial configuration
and $\xi_N^{(m)}(\cdot)=\sum_{i=1}^{N} 
\delta_{x_i^{(m)}}(\cdot) \in \mM^+_0$,
$\x^{(m)}=(x^{(m)}_1, \dots, x^{(m)}_N) \in \W_N^{+}$ for 
the configurations at times $t_m, 1 \leq m \leq M$.
In (\ref{eqn:pxi}), if some of $x_i$'s in 
$\x$ coincide, the factor
$f_N^{(\nu)}(t_1, \x^{(1)}|\x)/h_N(\x)$
is interpreted using l'H\^opital's rule.

For $\x^{(m)}=(x^{(m)}_1, \dots, x^{(m)}_N) \in \W_N^{+}$ with
$\xi_N^{(m)}(\cdot)=\sum_{i=1}^{N} \delta_{x^{(m)}_i} (\cdot)$
and $N' \in \{ 1,2, \dots, N\}$, we put
$\x^{(m)}_{N'}=(x^{(m)}_1, \dots, x^{(m)}_{N'})
\in \W_{N'}^{+}$,
$1 \leq m \leq M$.
For a sequence $(N_m)_{m=1}^{M}$ of positive integers 
less than or equal to $N$,
we define the 
$(N_1, \dots, N_{M})$-{\it multitime correlation function} by
\begin{eqnarray}
&& \rho^{\xi_N}_{\nu} \Big(t_{1}, \x^{(1)}_{N_1} ; 
\dots; t_M, \x^{(M)}_{N_M} \Big) 
\nonumber\\
&&=
\int_{\prod_{m=1}^{M} \R_+^{N-N_{m}}}
\prod_{m=1}^{M} \prod_{i=N_{m}+1}^{N} dx_{i}^{(m)}
p_{\nu}^{\xi_N} \Big(t_1, \xi_N^{(1)}; 
\dots; t_M, \xi_N^{(M)} \Big)
\prod_{m=1}^{M}
\frac{1}{(N-N_{m})!},
\label{eqn:corr}
\end{eqnarray}
which is symmetric in the sense that
$
\rho^{\xi_N}_{\nu}(\dots; t_m, \sigma(\x^{(m)}_{N_m}); \dots)
=\rho^{\xi_N}_{\nu}(\dots; t_m, \x^{(m)}_{N_m}; \dots)
$
with
$\sigma(\x^{(m)}_{N_m})
\equiv (x^{(m)}_{\sigma(1)}, \dots, x^{(m)}_{\sigma(N_m)})$
for any permutation $\sigma \in {\cal S}_{N_m},
1 \leq \forall m \leq M$. 
For any $M \in \N$,
$f_m \in {\rm C}_{0}(\R_+), \theta_m \in \R,
1 \leq m \leq M$,
$0 < t_1 < \cdots < t_M < \infty$,
the Laplace transform of (\ref{eqn:pxi}) is considered
as a functional of $\vchi(x)=(\chi_1(x), \dots, \chi_M(x))$,
where 
$\chi_m(x) \equiv e^{\theta_m f_m(x)}-1, 1 \leq m \leq M,
x \in \R_+$;
\begin{eqnarray}
\cG^{\xi_N}_{\nu}[\vchi] 
&\equiv& \int_{\R_+^{NM}} \prod_{m=1}^{M} d \x^{(m)} \,
p^{\xi_N}_{\nu} \Big( t_1, \xi_N^{(1)}; \dots ;
t_M, \xi_N^{(M)} \Big)
\exp \left\{ \sum_{m=1}^{M} \theta_m
\int_{\R} f_m(x) \xi_N^{(m)}(dx) \right\} 
\nonumber\\
&=& \E^{\xi_N} \left[
\exp \left\{ \sum_{m=1}^{M} \theta_m
\int_{\R} f_m(x) \Xi^{(\nu)}(t_m, dx) \right\} \right].
\label{eqn:G1}
\end{eqnarray}
It is the generating function of multitime correlation
functions, since if we expand it with respect to
$\chi_{m}(x_i^{(m)})$'s, 
$\rho_{\nu}^{\xi_N}$'s appear as coefficients in terms;
\begin{eqnarray}
\cG^{\xi_N}_{\nu}[\vchi]
&=& \sum_{N_{1}=0}^{N} \cdots
\sum_{N_{M}=0}^{N}
\prod_{m=1}^{M} \frac{1}{N_m !} 
\int_{\R_+^{N_m}} d \x^{(m)}_{N_m} \nonumber\\
&& \times 
\prod_{m=1}^{M} \prod_{i=1}^{N_{m}} 
\chi_{m}(x_{i}^{(m)})
\rho^{\xi_N}_{\nu} \Big( t_{1}, \x^{(1)}_{N_1};
\dots ; t_{M}, \x^{(M)}_{N_M} \Big).
\label{eqn:G2}
\end{eqnarray}

In the present paper, first we prove that, 
for any fixed initial configuration $\xi_N \in \mM^+$
with $\xi_N(\R_+) = N \in \N$, there is a function
$\mbK^{\xi_N}_{\nu}(s,x;t,y)$, which
is continuous with respect to
$(x,y) \in (0, \infty)^2$ for any fixed $(s,t) \in [0, \infty)^2$, 
and that the function (\ref{eqn:G1}) is given by
the {\it Fredholm determinant} in the form
\begin{equation}
\cG^{\xi_N}_{\nu}[\vchi]=
\Det_{
\substack{1 \leq m, n \leq M \\ (x,y) \in (0, \infty)^2}
}
 \Bigg[ \delta_{m n} \delta (x-y)
+\mbK^{\xi_N}_{\nu}(t_m,x;t_n,y)
\chi_n (y) \Bigg].
\label{eqn:Fredholm}
\end{equation}
By definition of Fredholm determinant 
(see Eq.(\ref{eqn:Fredholm2}) in Section 4.1),
(\ref{eqn:Fredholm}) means that any multitime correlation function
is given by a determinant
\begin{equation}
\rho^{\xi_N}_{\nu}
 \Big(t_1,\x^{(1)}_{N_1}; \dots;t_M,\x^{(M)}_{N_M} \Big) 
=\det_{
\substack{1 \leq i \leq N_{m}, 1 \leq j \leq N_{n} \\ 1 \leq m, n \leq M}
}
\Bigg[
\mbK^{\xi_N}_{\nu}(t_m, x_{i}^{(m)}; t_n, x_{j}^{(n)} )
\Bigg].
\label{eqn:rho1}
\end{equation}
The function $\mbK^{\xi_N}_{\nu}$ is called
the {\it correlation kernel} and it determines the finite
dimensional distributions of the process 
$(\Xi^{(\nu)}(t), \P^{\xi_N}_{\nu})$ 
through (\ref{eqn:rho1}).
It is an extension of {\it determinantal (Fermion)
point process} of distributions
studied by Soshnikov \cite{Sos00}
and Shirai and Takahashi \cite{ST03}
to the cases on ${\bf T} \times \R_+$
with ${\bf T}=\{t_1, \dots, t_M\}$,
$M \in \N, 0 < t_1< \dots < t_M < \infty$.
See \cite{HKPV09} for variety of examples 
of determinantal point processes. 
We express this result by simply saying
that the noncolliding BESQ$^{(\nu)}$
is {\it determinantal with a correlation kernel}
$\mbK^{\xi_N}_{\nu}$ for any $\xi_N \in \mM^+$ with
$\xi_N(\R_+)=N \in \N$ 
(Theorem \ref{thm:Finite}).

Next we consider the infinite-particle limits.
For $\xi \in \mM^+$ with $\xi(\R_+)= \infty$, when 
$\mbK_{\nu}^{\xi \cap [0, L]}$
converges to a continuous function as $L \to \infty$, 
the limit is written as $\mbK_{\nu}^{\xi}$. 
If
$
\sup_{x, y \in I}|
\mbK^{\xi \cap [0, L]}(s, x; t, y) | < \infty, \forall L >0
$
for any $(s,t) \in (0, \infty)^{2}$ and 
any compact interval $I \subset (0, \infty)$, 
we can obtain the convergence of generating functions
of multitime correlation functions,
$\cG_{\nu}^{\xi \cap [0, L]}[\vchi] \to
\cG_{\nu}^{\xi}[\vchi]$, as $L \to \infty$.
It implies 
$\P^{\xi \cap [0, L]}_{\nu} \to {^{\exists}}\P^{\xi}_{\nu}$ 
as $L \to \infty$
in the sense of finite dimensional distributions
weakly in the vague topology.
In this case, we say that the noncolliding BESQ$^{(\nu)}$
$(\Xi^{(\nu)}(t), \P^{\xi}_{\nu})$
with an infinite number of particles $\xi(\R_+)=\infty$
is {\it well defined with the correlation kernel}
$\mbK^{\xi}$ \cite{KT10}. 
We will give sufficient conditions so that
the process $(\Xi^{(\nu)}(t), \P^{\xi}_{\nu})$
is well defined, in which the correlation kernel
is generally expressed using a double
integral of an {\it entire function} represented by 
the Weierstrass canonical product
having zeros on $\supp \xi$,
where $\supp \xi = \{x \in \R : \xi(\{x\}) > 0\}$
(Theorem \ref{thm:Infinite}).
As an application of this theorem, we will study 
the following example of infinite particle system,
which is a non-equilibrium dynamics 
exhibiting a {\it relaxation phenomenon}.

Consider the Bessel function
\begin{equation}
J_{\nu}(z) = \sum_{n=0}^{\infty} 
\frac{(-1)^n}{\Gamma(n+1) \Gamma(n+1+\nu)}
\left( \frac{z}{2} \right)^{2n+\nu}.
\label{eqn:J1}
\end{equation}
It is an analytic function of $z$ in a cut plane.
The function $J_{\nu}(z)/z^{\nu}$ is an entire function.
As usual we define $z^{\nu}$
to be $\exp(\nu \log z)$, where the argument of $z$
is given its principal value;
\begin{equation}
z^{\nu}=\exp \Big[ \nu \Big\{
\log |z| + \sqrt{-1} {\rm arg} (z) \Big\} \Big], \quad
-\pi < {\rm arg} (z) \leq \pi.
\label{eqn:znu}
\end{equation}
The function $J_{\nu}(z)$ is analytically continued
outside this range of ${\rm arg}(z)$ so that
the relation
\begin{equation}
J_{\nu}(e^{m \pi \sqrt{-1}} z)
= e^{\nu m \pi \sqrt{-1}} J_{\nu}(z)
\label{eqn:Jcont}
\end{equation}
holds \cite{Wat44}.
If $\nu > -1$, $J_{\nu}(z)$ has an infinite number of pairs
of positive and negative zeros with the same 
absolute value,
which are all simple. We write the positive zeros
of $J_{\nu}(z)$ arranged in ascending order of 
the absolute value as
$$
0 < j_{\nu,1} < j_{\nu,2} < j_{\nu, 3} < \cdots.
$$
Explicitly $J_{\nu}(z)$ is expressed using the 
{\it infinite product
of the Weierstrass primary factors of genus zero} as
(see Chapter XV of \cite{Wat44}),
\begin{equation}
J_{\nu}(z)=\frac{(z/2)^{\nu}}{\Gamma(\nu+1)}
\prod_{i=1}^{\infty} \left(
1- \frac{z^2}{j_{\nu,i}^2} \right).
\label{eqn:entire1}
\end{equation}
The configuration in which every point of
the square of positive zero of $J_{\nu}(z)$
is occupied by one particle, denoted by
\begin{equation}
\xi_{J_{\nu}}^{\langle 2 \rangle}(\cdot)
=\sum_{i =1}^{\infty}
\delta_{j_{\nu, i}^2}(\cdot),
\label{eqn:xiJnu2}
\end{equation}
satisfies the conditions of Theorem \ref{thm:Infinite}.
We will determine the correlation kernel of
the noncolliding BESQ$^{(\nu)}$ starting from
$\xi_{J_{\nu}}^{\langle 2 \rangle}$ 
explicitly (Theorem \ref{thm:RelBessel} (i))
and prove that the process shows
a relaxation phenomenon to a stationary
process,
$$
(\Xi^{(\nu)}(t+\theta), \P^{\xi_{J_{\nu}}^{\langle 2 \rangle}}_{\nu})
\to (\Xi^{(\nu)}(t), \bP_{J_{\nu}})
\quad \mbox{as} \quad \theta \to \infty
$$
weakly in the sense of finite dimensional distributions
(Theorem \ref{thm:RelBessel} (ii)).
Here $(\Xi^{(\nu)}(t), \bP_{J_{\nu}})$ 
is the equilibrium dynamics, which is 
determinantal with the
correlation kernel
\begin{equation}
\bK_{J_{\nu}}(t-s,y|x) = \left\{
   \begin{array}{ll}
\displaystyle{
\int_{0}^{1} d u \,
e^{-2u(s-t)} J_{\nu}(2 \sqrt{u x})
J_{\nu}(2 \sqrt{u y})
} 
& \mbox{if} \quad s < t  \\
& \\
\displaystyle{
\frac{J_{\nu}(2 \sqrt{x}) \sqrt{y} J_{\nu}'(2 \sqrt{y})
-\sqrt{x} J_{\nu}'(2\sqrt{x}) J_{\nu}(2\sqrt{y})}{x-y}
}
& \mbox{if} \quad t=s \\
& \\
\displaystyle{
- \int_{1}^{\infty} d u \,
e^{-2u(s-t)} J_{\nu}(2 \sqrt{ux})
J_{\nu}(2 \sqrt{uy})
}
& \mbox{if} \quad s > t,
   \end{array} \right. 
\label{eqn:KBessel1}
\end{equation}
$(x, y) \in (0, \infty)^2$, 
where $J'(z)=dJ(z)/dz$.
Note that this kernel is temporally homogeneous but
spatially inhomogeneous.
Let $\mu_{J_{\nu}}$ be the determinantal (Fermion) point
process on $\R_+$, in which, for any $N \in \N$, 
$N$-point correlation function is given by
\cite{For10}
$$
\rho_{\nu}(\x_{N})=\det_{1 \leq i, j \leq N}
\Big[ K_{J_{\nu}}(x_i|x_j) \Big],
\quad \x_{N}=(x_1, \dots, x_{N}) \in (0, \infty)^N,
$$
with the {\it Bessel kernel}
\begin{eqnarray}
K_{J_{\nu}}(y|x) &\equiv& \bK_{J_{\nu}}(0, y|x)
\nonumber\\
&=& \frac{J_{\nu}(2 \sqrt{x}) \sqrt{y} J_{\nu}'(2 \sqrt{y})
-\sqrt{x} J_{\nu}'(2\sqrt{x}) J_{\nu}(2\sqrt{y})}{x-y}
\nonumber\\
&=& 
\frac{\sqrt{x} J_{\nu+1}(2 \sqrt{x}) J_{\nu}(2 \sqrt{y})
- J_{\nu}(2 \sqrt{x}) \sqrt{y} J_{\nu+1}(2 \sqrt{y})}
{x-y}.
\label{eqn:KJ0}
\end{eqnarray}
In particular, the density of particle
at $x \in (0, \infty)$ is given by
\begin{eqnarray}
\rho_{\nu}(x) &=& K_{J_{\nu}}(x|x)
\equiv \lim_{y \to x} K_{J_{\nu}}(y|x)
\nonumber\\
&=& (J_{\nu}'(2 \sqrt{x}))^2
+\left( 1- \frac{\nu^2}{4x} \right)
(J_{\nu}(2 \sqrt{x}))^2.
\nonumber\\
&=& (J_{\nu}(2 \sqrt{x}))^2
-J_{\nu+1}(2 \sqrt{x}) J_{\nu-1}(2 \sqrt{x}).
\label{eqn:rho0}
\end{eqnarray}
The probability measure $\mu_{J_{\nu}}$ is obtained
in an $N \to \infty$ limit called the {\it hard-edge scaling-limit}
of the distribution of squares of eigenvalues
of random matrices in the chiral Gaussian unitary 
ensemble \cite{Meh04,For10}. 
The process $(\Xi^{(\nu)}(t), \bP_{J_{\nu}})$ is 
a reversible process with respect to $\mu_{J_{\nu}}$.
The correlation kernel (\ref{eqn:KBessel1}) is
called the {\it extended Bessel kernel} \cite{FNH99,TW04}.

In the random matrix theory, three kinds of determinantal
point processes of infinite particle systems
have been well studied, the correlation kernels
of which are given by (i) the sine kernel,
(ii) the Airy kernel, and (iii) the Bessel kernel
\cite{Meh04,For10}.
They are obtained by taking 
(i) the bulk, (ii) the soft-edge, and (iii)
the hard-edge scaling limits in the
Gaussian unitary ensemble (GUE) for (i) and (ii)
and in the chiral GUE for (iii), respectively.
These three determinantal point processes
have been extended to time-dependent versions
so that they describe equilibrium dynamics,
which are reversible with respect to the determinantal 
point processes \cite{KT07b}.
Corresponding to these three stationary processes,
the present authors introduced three relaxation processes
with infinite numbers of particles realized in
(i) the Dyson model (the noncolliding Brownian motion)
starting from $\Z$ ({\it i.e.} the zeros of
$\sin (\pi z)$) in \cite{KT10},
(ii) the Dyson model with drift terms
starting from the Airy zeros in \cite{KT09},
and (iii) the noncolliding BESQ$^{(\nu)}$ starting from
the squares of positive zeros of $J_{\nu}$ in the present paper.
The scaling limits are performed by increasing
the number of particles in the system $N \to \infty$,
while in our setting determinantal processes with
infinite numbers of particles $N=\infty$ are
well constructed based on the theory of entire functions.
In our relaxation processes in non-equilibrium
we only have to wait for sufficiently long time
to observe the three determinantal point processes.
In order to give temporally inhomogeneous correlation
kernels explicitly, we have reported the relaxation
processes with the special initial configurations.
The universality of the three determinantal point
processes in a wide variety of
fields of mathematics, physics, and others \cite{Meh04,For10}
implies robustness of relaxation phenomena
with respect to initial configurations.
Mathematical justification of this fact 
will be reported in the future.

The present paper is organized as follows.
In Section 2 preliminaries and main results
are given.
In Section 3 the properties of special functions
used in this paper are given.
Section 4 is devoted to proofs of results.

\SSC{Preliminaries and Main Results}

We introduce the following operations; 
for $\xi(\cdot)=\sum_{i\in \Lambda}
\delta_{x_i}(\cdot) \in\mM$,
\begin{description}
\item[(shift)] with $u \in \R$, 
$\tau_u \xi(\cdot) =\displaystyle{\sum_{i \in \Lambda}} 
\delta_{x_i+u}(\cdot)$,

\item[(dilatation)] with $c>0$,
$c \circ \xi(\cdot)=\displaystyle{\sum_{i \in \Lambda} 
\delta_{c x_i}(\cdot)}$,

\item[(square)]
$\displaystyle{
\xi^{\langle 2 \rangle}(\cdot)
=\sum_{i \in \Lambda} \delta_{x_i^2} (\cdot)}$, 
\end{description}
and for $\xi(\cdot)=\sum_{i\in \Lambda}
\delta_{x_i}(\cdot) \in\mM^+$,
\begin{description}
\item[(square root)]
$\displaystyle{
\xi^{\langle 1/2 \rangle}(\cdot)
=\sum_{i \in \Lambda} 
\Big( \delta_{\sqrt{x_i}}+ \delta_{-\sqrt{x_i}} (\cdot) \Big)
}$.
\end{description}

Note that the notation (\ref{eqn:xiJnu2}) states that
this configuration is obtained as the square
of the point-mass distribution on the positive
zeros of $J_{\nu}(z)$ denoted by
$\xi_{J_{\nu}}(\cdot)=\sum_{i=1}^{\infty} 
\delta_{j_{\nu,i}}(\cdot)$.
We use the convention such that
$$
\prod_{x\in\xi}f(x) =\exp
\left\{\int_\R \xi(dx) \log f(x) \right\}
=\prod_{x \in \supp \xi}f(x)^{\xi(\{x\})}
$$
for $\xi\in \mM$ and a function $f$ on $\R$.
For a multivariate symmetric function $g$ we write 
$g((x)_{x \in \xi})$ for $g((x_i)_{i \in \Lambda})$.

The transition probability density of
BESQ$^{(\nu)}$, given by (\ref{eqn:p0}), satisfies the 
Chapman-Kolmogorov equation
\begin{equation}
\int_{0}^{\infty} dy \,
p^{(\nu)}(t-s, z|y) p^{(\nu)}(s, y|x)
= p^{(\nu)}(t,z|x),
\label{eqn:CK1}
\end{equation}
for $0 \leq s \leq t, x, z \in \R_+$.
Here we define the modified Bessel function
of the first kind on $\C$ as
\begin{equation}
I_{\nu}(z) = \left\{ 
\begin{array}{ll}
e^{-\nu \pi \sqrt{-1}/2} 
J_{\nu}(e^{\pi \sqrt{-1}/2} z),
& - \pi < {\rm arg} (z) \leq \pi/2, \cr
e^{3 \nu \pi \sqrt{-1}/2} 
J_{\nu}(e^{-3 \pi \sqrt{-1}/2} z),
& \pi/2 < {\rm arg} (z) \leq \pi,
\end{array} \right.
\label{eqn:Icont}
\end{equation}
where $J_{\nu}$ is defined by (\ref{eqn:J1})
so that (\ref{eqn:Jcont}) holds \cite{AAR99}.
This definition is consistent with (\ref{eqn:I1}),
associated with the relation
\begin{equation}
I_{\nu}(e^{m \pi \sqrt{-1}} z)
= e^{\nu m \pi \sqrt{-1}} I_{\nu}(z).
\label{eqn:Icont2}
\end{equation}
For $t \in \R$, we define
\begin{equation}
p^{(\nu)}(t, y|x)
=  \left\{ \begin{array}{ll}
\displaystyle{
\frac{1}{2|t|} \left( \frac{y}{x} \right)^{\nu/2}
\exp \left( - \frac{x+y}{2t} \right)
I_{\nu} \left( \frac{\sqrt{xy}}{t} \right)},
& \quad t \in \R \setminus \{0\}, x \in \C \setminus \{0\}, \cr
\displaystyle{
\frac{y^{\nu}}{(2 |t|)^{\nu+1} \Gamma(\nu+1)} e^{-y/2t}},
& \quad t \in \R \setminus \{0\}, x=0, \cr
& \cr
\delta(y-x),
& \quad t=0, x \in \C,
\end{array} \right.
\label{eqn:pnu-}
\end{equation}
$y \in \C$.
Then by using the modified version of Weber's integral
of the Bessel functions (see Eqs.(\ref{eqn:Watson2})
and (\ref{eqn:Int2}) in Section 3.1),
(\ref{eqn:CK1}) is extended to the
following equations.
For $0 \leq s \leq t, x, z \in \C$,
\begin{eqnarray}
\label{eqn:pnu-2}
&& \int_0^{\infty} dy \, 
p^{(\nu)}(-t, z|y) p^{(\nu)}(t-s, y|x)
=p^{(\nu)}(-s, z|x), \\
\label{eqn:pnu-2b}
&& \int_0^{\infty} dy \, 
p^{(\nu)}(t-s, z|y)
p^{(\nu)}(-t, y|x) 
=p^{(\nu)}(-s, z|x).
\end{eqnarray}
We define
\begin{equation}
p_{J_{\nu}}(t, y|x)=
\left\{ \begin{array}{ll}
\displaystyle{
\left( \frac{x}{y} \right)^{\nu/2}
p^{(\nu)}(t, y|x) }, 
& \quad t \in \R, x \in \C \setminus \{0\},
\cr
& \cr
y^{-\nu/2} p^{(\nu)}(t, y|0),
& \quad t \in \R, x=0,
\end{array} \right.
\label{eqn:pnupJ}
\end{equation}
$y \in \C$.
When $t \geq 0, x, y \in \R_+$, it has the expression
\begin{eqnarray}
p_{J_{\nu}}(t, y|x) 
&=& \int_{0}^{\infty} du \,
J_{\nu}(2 \sqrt{ux}) J_{\nu}(2 \sqrt{uy}) e^{-2 u t}
\nonumber\\
&=& 2 \int_{0}^{\infty} dw \, w
J_{\nu}(2 w \sqrt{x}) J_{\nu}(2 w \sqrt{y}) e^{-2 w^2 t}.
\label{eqn:pJnu1}
\end{eqnarray}
The Chapman-Kolmogorov equation (\ref{eqn:CK1}) 
and its extensions (\ref{eqn:pnu-2}) and (\ref{eqn:pnu-2b})
are mapped to
\begin{eqnarray}
\label{eqn:CKJ1}
&& \int_{0}^{\infty} dy \,
p_{J_{\nu}}(t-s, z|y) p_{J_{\nu}}(s,y|x)
=p_{J_{\nu}}(t, z|x), \\
\label{eqn:CKJ2}
&& \int_{0}^{\infty} dy \,
p_{J_{\nu}}(-t, z|y) p_{J_{\nu}}(t-s,y|x)
=p_{J_{\nu}}(-s, z|x), \\
\label{eqn:CKJ3}
&& \int_{0}^{\infty} dy \,
p_{J_{\nu}}(t-s, z|y) p_{J_{\nu}}(-t,y|x)
=p_{J_{\nu}}(-s, z|x)
\end{eqnarray}
for $0 \leq s \leq t, x, z \in \C$.

For $\xi_N \in \mM^+$ with $\xi_N(\R_+)=N \in \N$,
we define the functions of $z \in \C$,
\begin{eqnarray}
\label{eqn:entire2b}
\Pi_0(\xi_N,z) &=& 
\prod_{x \in \xi_N \cap \{0\}^{\rm c}}\left(1-\frac{z}{x} \right), 
\\
\label{eqn:entire2}
\Pi^{(\nu)}(\xi_N, z)
&=& z^{\nu/2} \Pi_0(\xi_N,z).
\end{eqnarray}
For $a \in \C$ we also define
\begin{eqnarray}
\label{eqn:entire3a}
\Phi_0(\xi_N, a,z) &=& \Pi_0(\tau_{-a} \xi_N, z-a)
= \prod_{x \in \xi_N \cap \{ a \}^{\rm c}}
\left(1-\frac{z-a}{x-a} \right), 
\\
\label{eqn:entire3}
\Phi^{(\nu)}(\xi_N,a, z)
&=& \left\{ \begin{array}{ll}
\displaystyle{\left( \frac{z}{a} \right)^{\nu/2}
\Phi_0(\xi_N, a,z)},
& \quad \mbox{if $a \not=0$}, \cr
& \cr
\Pi^{(\nu)}(\xi_N, z),
& \quad \mbox{if $a=0$}.
\end{array} \right.
\end{eqnarray}

\begin{thm}
\label{thm:Finite}
{\rm (i)} For any fixed configuration $\xi_N \in \mM^+$ 
with $\xi_N(\R_+)=N \in \N$,
$(\Xi^{(\nu)}(t), \P_{\nu}^{\xi_N})$
is determinantal with the correlation kernel
\begin{eqnarray}
\mbK^{\xi_N}_{\nu}(s, x; t, y)
&=& \lim_{\varepsilon \downarrow 0}
\frac{1}{2 \pi \sqrt{-1}} 
\int_{-\infty}^{-\varepsilon} dy' 
\oint_{\Gamma_{y'}(\xi_N)} dz \,  
p^{(\nu)}(s, x|z) \frac{1}{y'-z}
\Phi_0(\xi_N, z, y')
p^{(\nu)}(-t, y'|y)
\nonumber\\
&& - {\bf 1}(s > t)p^{(\nu)}(s-t, x|y),
\quad (s,t) \in [0, \infty)^2, (x,y) \in (0, \infty)^2, 
\label{eqn:KN1a}
\end{eqnarray}
where $\Gamma_{y'}(\xi_N)$ denotes a 
counterclockwise contour on the
complex plane $\C$ encircling the points 
in $\supp \xi_N$ on $\R_+$
but not the point $y' \in (-\infty, -\varepsilon]$, 
and ${\bf 1}(\omega)$ is the indicator function
of condition $\omega$. \\
{\rm (ii)} 
If $\xi_N \in \mM^+_0$ with
$\xi_N(\R_+)=N \in \N$, the correlation kernel is
given by
\begin{eqnarray}
\mbK^{\xi_N}_{\nu}(s, x; t, y)
&=&  
\int_{0}^{\infty} \xi_N(dx') 
\int_{-\infty}^{0} dy' \,
 p^{(\nu)}(s, x|x') \Phi_0(\xi_N, x', y')
 p^{(\nu)}(-t, y'|y)
\nonumber\\
&& - {\bf 1}(s > t) p^{(\nu)}(s-t, x|y),
\quad (s,t) \in [0, \infty)^2, (x,y) \in (0, \infty)^2.
\label{eqn:KN1bz}
\end{eqnarray}
Without changing any finite dimensional distributions
of the process, 
the correlation kernel (\ref{eqn:KN1bz}) can be
replaced by 
\begin{eqnarray}
\mbK^{\xi_N}_{J_{\nu}}(s, x; t, y)
&=&  
\int_{0}^{\infty} \xi_N(dx') 
\int_{-\infty}^{0} dy' \,
 p_{J_{\nu}}(s, x|x') \Phi^{(\nu)}(\xi_N, x', y')
  p_{J_{\nu}}(-t, y'|y)
\nonumber\\
&& - {\bf 1}(s > t)p_{J_{\nu}}(s-t, x|y),
\quad (s,t) \in [0, \infty)^2, (x,y) \in (0, \infty)^2.
\label{eqn:KN1b}
\end{eqnarray}
\label{eqn:THfinite}
\end{thm}
\noindent{\bf Remark 1.} 
If we consider a spatial distribution of particles
at a single time $t \geq 0$ of the noncolliding BESQ$^{(\nu)}$,
we have a determinantal point process
with the correlation function $\mbK^{\xi_N}_{\nu}(t,x; t,y),
(x,y) \in (0, \infty)^2$.
In particular, if we set $s=t=1/2$,
the correlation kernel (\ref{eqn:KN1a}) is 
reduced to be the kernel of 
the perturbed chiral GUE
of random matrices given in Proposition 5
by Desrosiers and Forrester \cite{DF08}.
More detail, see Remark 2 in Section 3.3. 
\vskip 0.5cm

For $L > 0, \alpha >0$ and $\xi \in \mM$ we put
$$
M_{\alpha}(\xi, L)=\left( \int_{[-L, L] \setminus \{0\}}
\frac{\xi(dx)}{|x|^{\alpha}} \right)^{1/\alpha}
$$
and
$$
M_{\alpha}(\xi)=\lim_{L \to \infty} 
M_{\alpha}(\xi, L),
$$
if the limit finitely exists.
We introduce the following conditions
for configurations $\xi \in \mM^+$.
\vskip 0.3cm
\noindent
{\bf (C.A)} (i)
There exists $\alpha \in (1/2,1)$ and $C_1 > 0$
such that 
$M_{\alpha}(\xi) \leq C_1$. 

(ii) There exist $\beta > 0$ and $C_2 > 0$
such that
$$
M_1(\tau_{-a} \xi) 
\leq C_2 (|a| \vee 1)^{-\beta},
\quad \forall a \in \supp \xi.
$$
\vskip 0.3cm
\noindent
We denote by $\mX^+$ the set of configurations
satisfying the conditions {\bf (C.A)}, and put
$\mX_0^+=\mX^+ \cap \mM_0^+$.
For $\xi \in \mX_0^+, a \in \R$ and $z \in \C$ we can define
\begin{eqnarray}
\Phi_0(\xi, a, z)
&=& \lim_{L \to \infty} \Phi_0(\xi \cap [a-L, a+L], a, z),
\nonumber\\
\Phi^{(\nu)}(\xi, a, z)
&=& \lim_{L \to \infty} \Phi^{(\nu)}(\xi \cap [a-L, a+L], a, z),
\quad \nu > -1.
\nonumber
\end{eqnarray}
Since $\Phi_0(\xi, a, z)$ has the expression of
the {\it Weierstrass canonical product of genus zero}
(see (\ref{eqn:entire3a})),
it is an entire function of a variable $z \in \C$ \cite{Lev96}.
The set of zeros of $\Phi^{(\nu)}(\xi, a, z)$
is given by $(\supp \xi \cup \{0\}) \cap \{a\}^{\rm c}$
and all zeros except 0 are simple for $\xi \in \mM_0^+$.

\begin{thm}
\label{thm:Infinite}
If $\xi \in \mX_0^+$, the process
$(\Xi^{(\nu)}(t), \P_{\nu}^{\xi})$ is well defined
with the correlation kernel 
\begin{eqnarray}
\mbK^{\xi}_{\nu}(s, x; t, y)
&=&  
\int_{0}^{\infty} \xi(dx') 
\int_{-\infty}^{0} dy' \,
 p^{(\nu)}(s, x|x') \Phi_0(\xi, x', y')
 p^{(\nu)}(-t, y'|y)
\nonumber\\
&& - {\bf 1}(s > t) p^{(\nu)}(s-t, x|y),
\quad (s,t) \in [0, \infty)^2, (x,y) \in (0, \infty)^2.
\label{eqn:K1bz}
\end{eqnarray}
This correlation kernel $\mbK^{\xi}_{\nu}$
can be replaced by
\begin{eqnarray}
\mbK^{\xi}_{J_{\nu}}(s, x; t, y)
&=&  
\int_{0}^{\infty} \xi(dx') 
\int_{-\infty}^{0} dy' \,
 p_{J_{\nu}}(s, x|x') \Phi^{(\nu)}(\xi, x', y')
  p_{J_{\nu}}(-t, y'|y)
\nonumber\\
&& - {\bf 1}(s > t)p_{J_{\nu}}(s-t, x|y),
\quad (s,t) \in [0, \infty)^2, (x,y) \in (0, \infty)^2
\label{eqn:K1b}
\end{eqnarray}
without changing any finite dimensional distributions
of the process.
\label{eqn:THinfinite}
\end{thm}
In case $\xi(\R_+)=\infty$, Theorem \ref{thm:Infinite}
gives the noncolliding BESQ$^{(\nu)}$ 
with an infinite number of particles
starting from the configuration $\xi \in \mX^+_0$.
It is easy to check that,
if $(x,y) \in (0, \infty)^2$,
$$
\mbK^{\xi}_{\nu}(t, x; t, y) \mbK^{\xi}_{\nu}(t, y; t, x) dxdy 
\to \xi(dx){\bf 1}(x=y), 
\quad t \to 0 \quad
\mbox{in the vague topology}.
$$
For $\gamma >0$, we put
$$
g^{\gamma}(x)=x^{\gamma}, \quad x \in \R_+,
$$
and
$$
  \eta^{\gamma}(\cdot)=\sum_{i=1}^{\infty}
  \delta_{g^{\gamma}(i)} (\cdot).
$$
For any $\gamma > 1$ we can show by simple calculation
that $\eta^{\gamma}$ satisfies {\bf (C.A)}(i) with
any $\alpha \in (1/\gamma, 1)$
and some $C_1=C_1(\alpha) > 0$ depending on $\alpha$,
and {\bf (C.A)} (ii) with any $\beta \in (0, \gamma-1)$
and some $C_2=C_2(\beta) > 0$ depending on $\beta$.
This implies that $\eta^{\gamma}$ is an element
of $\mX^+_0$ for any $\gamma > 1$.

More interesting example is given by the following
theorem.

\begin{thm}
\label{thm:RelBessel}
{\rm (i)} 
The noncolliding BESQ$^{(\nu)}$ starting from
$\xi_{J_{\nu}}^{\langle 2 \rangle}$,
$(\Xi^{(\nu)}(t), \P^{\xi_{J_{\nu}}^{\langle 2 \rangle}}_{\nu})$,
is well defined with the correlation kernel
\begin{eqnarray}
\mbK_{J_{\nu}}(s,x;t,y)
&=& \sum_{i=1}^{\infty}
\int_{-\infty}^0 dz \,
p_{J_{\nu}}(s, x|j_{\nu, i}^2)
\frac{2j_{\nu, i}}{j_{\nu, i}^2-z}
\frac{J_{\nu}(\sqrt{z})}{J_{\nu+1}(j_{\nu, i})}
p_{J_{\nu}}(-t, z|y)
\nonumber\\
&& - {\bf 1}(s>t)
p_{J_{\nu}}(s-t, x|y),
\quad (s,t) \in [0, \infty)^2,
(x,y) \in (0, \infty)^2.
\label{eqn:KJnu}
\end{eqnarray}
{\rm (ii)} Let $(\Xi^{(\nu)}(t), \bP_{J_{\nu}})$ be the 
equilibrium dynamics, which is determinantal with
the extended Bessel kernel (\ref{eqn:KBessel1}).
Then, for $t \geq 0$
\begin{equation}
(\Xi^{(\nu)}(t+ \theta), \P^{\xi_{J_{\nu}}^{\langle 2 \rangle}}_{\nu})
\to (\Xi^{(\nu)}(t), \bP_{J_{\nu}})
\quad \mbox{as} \quad \theta \to \infty
\label{eqn:relaxBessel}
\end{equation}
weakly in the sense of finite dimensional 
distributions.
\end{thm}

\SSC{Some Properties of Special Functions}
\subsection{Integral formulas of Bessel functions}

The following integral formulas are
known \cite{Wat44,AAR99}.
For $\Re \nu >-1, p, a, b > 0$,
\begin{eqnarray}
\label{eqn:Watson2}
&&\int_{0}^{\infty} du \, u e^{-p^2 u^2}
J_{\nu}(au) J_{\nu}(bu)
=\frac{1}{2 p^2} 
\exp \left( -\frac{a^2+b^2}{4 p^2} \right)
I_{\nu} \left( \frac{ab}{2 p^2} \right), \\
\label{eqn:Int2}
&& \int_{0}^{\infty} du \, u e^{-p^2 u^2}
I_{\nu}(au) I_{\nu}(bu)
=\frac{1}{2 p^2} 
\exp \left( \frac{a^2+b^2}{4 p^2} \right)
I_{\nu} \left( \frac{ab}{2 p^2} \right).
\end{eqnarray}
The equalities (\ref{eqn:CK1}), (\ref{eqn:pnu-2}) and 
(\ref{eqn:pnu-2b})
for $p^{(\nu)}$
and (\ref{eqn:pJnu1})-(\ref{eqn:CKJ3}) 
for $p_{J_{\nu}}$ 
are derived from the above integral formulas.
In addition to them, the following 
equality is also derived from (\ref{eqn:Watson2})
with (\ref{eqn:Icont}) and (\ref{eqn:pJnu1}).
For $t > 0, x, z >0$
\begin{equation}
\int_{-\infty}^{0} dy \,
J_{\nu}(\sqrt{zy}) p_{J_{\nu}}(-t, y|x)
=e^{tz/2} J_{\nu}(\sqrt{zx}).
\label{eqn:pJnu-}
\end{equation}

\begin{lem}
\label{thm:j1}
For any $i \in \N, z \not= j_{\nu, i}$
\begin{eqnarray}
\label{eqn:j1}
\frac{j_{\nu,i}}{j_{\nu, i}^2-z^2}
\frac{J_{\nu}(z)}{J_{\nu+1}(j_{\nu, i})}
&=&\frac{1}{(J_{\nu+1}(j_{\nu, i}))^2}
\int_{0}^{1} du \, u J_{\nu}(zu) J_{\nu}(j_{\nu, i} u)
\\
\label{eqn:j2}
&=& \frac{1}{(J_{\nu+1}(j_{\nu, i}))^2} \frac{1}{2}
\int_{0}^{1} dw \, J_{\nu}(z \sqrt{w}) 
J_{\nu}(j_{\nu, i} \sqrt{w}).
\end{eqnarray}
\end{lem}
\noindent{\it Proof.} \quad
The following formula is found on page 482 in \cite{Wat44}
$$
\int_0^x t J_{\nu}(\alpha t) J_{\nu}(\alpha_0 t) dt
=\frac{x}{\alpha^2-\alpha_0^2}
\Big[ J_{\nu}(\alpha x) \alpha_0 J'_{\nu}(\alpha_0 x)
-J_{\nu}(\alpha_0 x) \alpha J'_{\nu}(\alpha x) \Big].
$$
Set
$\alpha_0=j_{\nu, i}$ and $x=1$.
Since $J_{\nu}(j_{\nu, i})=0$ by definition
of $j_{\nu, i}$'s, we have
$$
\int_{0}^{1} t J_{\nu}(\alpha t) J_{\nu}(j_{\nu, i} t) dt
= \frac{j_{\nu,i}}{\alpha^2-j_{\nu, i}^2}
J_{\nu}(\alpha) J'_{\nu}(j_{\nu, i}).
$$
Note that the Bessel function satisfies the relation
$
J'_{\nu}(z)=(\nu/z) J_{\nu}(z)- J_{\nu+1}(z).
$
Then
\begin{equation}
J'_{\nu}(j_{\nu, i})=-J_{\nu+1}(j_{\nu, i})
\label{eqn:J'2}
\end{equation}
and thus
$$
\int_{0}^{1} t J_{\nu}(\alpha t) J_{\nu}(j_{\nu, i} t) dt
= \frac{j_{\nu,i}}{j_{\nu, i}^2 - \alpha^2}
J_{\nu}(\alpha) J_{\nu+1}(j_{\nu, i}).
$$
Divide the both sides by $(J_{\nu+1}(j_{\nu, i}))^2$,
we obtain (\ref{eqn:j1}).
Eq. (\ref{eqn:j2}) is obtained from (\ref{eqn:j1})
by setting $u=\sqrt{w}$. \qed
\vskip 0.3cm

\subsection{Fourier-Bessel expansion}

As shown in Section 18.24 in \cite{Wat44},
any continuous function $f(x) \in L^2(0,1)$ has the expansion
$$
f(x)=\sum_{i=1}^{\infty} a_{i} 
J_{\nu}(j_{\nu, i} x)
\quad \mbox{with} \quad
a_{i}=\frac{2}{(J_{\nu+1}(j_{\nu, i}))^2}
\int_0^{1} u f(u) J_{\nu}(j_{\nu, i} u) du.
$$
That is,
$$
f(x)=\int_{0}^{1} du \, u f(u)
\sum_{i=1}^{\infty} \frac{2 J_{\nu}(j_{\nu, i} x)
J_{\nu}(j_{\nu, i} u)}
{(J_{\nu+1}(j_{\nu, i}))^2},
\quad f \in L^2(0, 1),
$$
which is called the {\it Fourier-Bessel expansion}.
Set $u=\sqrt{y}$ and then replace $x$ by $\sqrt{x}$,
we have
$$
f(\sqrt{x})=\int_{0}^{1} dy \, f(\sqrt{y})
\sum_{i=1}^{\infty} \frac{J_{\nu}(j_{\nu, i} \sqrt{x})
J_{\nu}(j_{\nu, i}\sqrt{y})}
{(J_{\nu+1}(j_{\nu, i}))^2}.
$$
In other words, the functions
$\{J_{\nu}(j_{\nu, i} \sqrt{x})/J_{\nu+1}(j_{\nu,i}), i \in \N\}$
form an orthonormal basis for $f \in L^2(0,1)$ 
and the completeness is also established;
\begin{equation}
\sum_{i=1}^{\infty} \frac{J_{\nu}(j_{\nu, i} \sqrt{x})
J_{\nu}(j_{\nu, i} \sqrt{y})}
{(J_{\nu+1}(j_{\nu, i}))^2}
=\delta(x-y),
\quad x, y \in (0,1).
\label{eqn:Fourier-Bessel}
\end{equation}

\subsection{Multiple orthogonal polynomials}

Fix a configuration
$$
\xi_N \in \mM^+ \quad
\mbox{with} \quad
\xi_N(\R_+)=N \in \N.
$$
In the present paper
the multiple orthogonal polynomials associated with
the modified Bessel function $I_{\nu}$ 
indexed by $\xi_N$ are
defined by the following \cite{CA03a,CA03b,KMW09,DF08}.

\begin{description}
\item[Type I:] \,
The {\it multiple orthogonal polynomials of the type I}
are the set of functions
\begin{equation}
\Big\{ A_{\xi_N}^{(\nu)}(y,x) : x \in \supp \xi_N,
\mbox{polynomial of $y$ of degree $\xi_N(x)-1$} \Big\}
\label{eqn:type1a}
\end{equation}
such that, if we set
\begin{equation}
Q_{\xi_N}^{(\nu)}(y)= \sum_{x \in \supp \xi_N}
A_{\xi_N}^{(\nu)}(y,x) 
\frac{1}{2} \left( \frac{y}{x} \right)^{\nu/2}
I_{\nu}(\sqrt{xy}) e^{-(x+y)/2},
\label{eqn:type1b}
\end{equation}
then
\begin{equation}
\int_{0}^{\infty} dy \, Q_{\xi_N}^{(\nu)}(y) y^{i}
= \left\{
\begin{array}{ll}
0, & i=0,1, \dots, \xi_N(\R_+)-2 \cr
1, & i=\xi_N(\R_+)-1.
\end{array} \right.
\label{eqn:type1c}
\end{equation}

\item[Type II:] \,
The {\it multiple orthogonal polynomial of the type II}
is the monic polynomial of degree $\xi_N(\R_+)$,
\begin{equation}
P_{\xi_N}^{(\nu)}(y)=y^{\xi_N(\R_+)} +
{\cal O}(y^{\xi_N(\R_+)-1})
\label{eqn:type2a}
\end{equation}
such that for each $x \in \supp \xi_N$
\begin{equation}
\int_{0}^{\infty} dy \,
P_{\xi_N}^{(\nu)}(y) y^{i}
\frac{1}{2} \left( \frac{y}{x} \right)^{\nu/2}
I_{\nu}(\sqrt{xy}) e^{-(x+y)/2}=0,
\quad 0 \leq i \leq \xi_N(x)-1.
\label{eqn:type2b}
\end{equation}
\end{description}

The following integral representations 
have been obtained by Desrosiers and Forrester \cite{DF08}.
\begin{lem}
\label{thm:QPintegral}
The functions $Q_{\xi_N}^{(\nu)}(y)$
and $P_{\xi_N}^{(\nu)}(y)$ have the 
following integral representations,
\begin{eqnarray}
\label{eqn:Qint}
&& Q_{\xi_N}^{(\nu)}(y)=
\frac{1}{2 \pi \sqrt{-1}}
\oint_{\Gamma(\xi_N)} dz \,
\frac{1}{2} \left( \frac{y}{z} \right)^{\nu/2}
I_{\nu}(\sqrt{yz}) e^{-(y+z)/2}
\frac{1}{\prod_{x \in \xi_N}(z-x)}, \\
\label{eqn:Pint}
&& P_{\xi_N}^{(\nu)}(y)=
\int_{-\infty}^{0} dw \,
\frac{1}{2} \left( \frac{w}{y} \right)^{\nu/2}
I_{\nu}(-\sqrt{yw}) e^{(y+w)/2}
\prod_{x \in \xi_N}(w-x),
\end{eqnarray}
where $\Gamma(\xi_N)$ denotes a 
counterclockwise contour on the
complex plane $\C$ encircling the points 
in $\supp \xi_N$ on $\R_+$.
\end{lem}
\noindent{\bf Remark 2.} 
The present definition of multiple orthogonal 
polynomials associated with
the modified Bessel function $I_{\nu}$
is slightly different from that
given by Desrosiers and Forrester \cite{DF08}.
Moreover, since
they have used the function
$
{_0}F{_1}(\alpha; z)
=\sum_{n=0}^{\infty} z^n/\{(\alpha)_n n!\}
$
in order to express the polynomials
instead of $I_{\nu}(z)$,
our expressions (\ref{eqn:Qint}) and (\ref{eqn:Pint})
seem to be quite different from their functions.
The identity
$$
I_{\nu}(z)=\frac{(z/2)^{\nu}}{\Gamma(\nu+1)}
{_0}F{_1} \left( \nu+1; \frac{z^2}{4} \right)
$$
is established, however, 
and then we can see that Lemma \ref{thm:QPintegral}
given above is equivalent with 
Proposition 6 of \cite{DF08}.
More precisely speaking, if we write the
orthogonal polynomials in \cite{DF08} as
$Q_{\vec{n}}^{\rm DF}(y)$ and $P_{\vec{n}}^{\rm DF}(y)$,
where the parameters $\alpha=\nu$, 
$\vec{n}=(n_1, \dots, n_D)$, and
${\bf a}={\bf b}^{\vec{n}}$,
we have the correspondence
\begin{equation}
Q_{\xi_N}^{(\nu)}(y)= 2^{-|\vec{n}|}
Q_{\vec{n}}^{\rm DF} \left( \frac{y}{2} \right),
\qquad
P_{\xi_N}^{(\nu)}(y) = 2^{|\vec{n}|}
P_{\vec{n}}^{\rm DF} \left( \frac{y}{2} \right),
\label{eqn:corresp}
\end{equation}
for $\xi_N(\cdot)=\sum_{i=1}^{N} \delta_{x_i}(\cdot)
=\sum_{\ell=1}^{D} n_{\ell} \delta_{2 b_{\ell}}(\cdot)$ and
$\xi_N(\R_+)=N=\sum_{\ell=1}^D n_{\ell}=|\vec{n}|$.
\vskip 0.5cm

We write $\xi_N(\cdot) = \sum_{i=1}^{N} \delta_{x_i}(\cdot)$
with $\x=(x_i)_{i=1}^{N}$ such that
$0 \leq x_1 \leq x_2 \leq \cdots \leq x_N$.
Then we set
$$
\xi_{N, 0}(\cdot) \equiv 0 \quad
\mbox{and} \quad
\xi_{N,i}(\cdot)= \sum_{k=1}^{i} \delta_{x_k}(\cdot),
\quad 1 \leq i \leq  N.
$$
By this definition
$\xi_{N, i}(\R_+)=i, 0 \leq i \leq N$
and $\xi_{N,i}(\{x\}) \leq \xi_{N, i+1}(\{x\}),
\forall x \in \R_+, 0 \leq i \leq N-1$. 
We define
\begin{equation}
M^{(\nu, +)}_i(y; \xi_N)=Q^{(\nu)}_{\xi_{N,i+1}}(y), 
\quad
M^{(\nu, -)}_i(y; \xi_N)=P^{(\nu)}_{\xi_{N,i}}(y),
\quad 0 \leq i \leq N-1.
\label{eqn:M+-def}
\end{equation}
By the orthogonality relations (\ref{eqn:type1c}),
(\ref{eqn:type2b}) and the above definitions, we can prove
the {\it biorthonormality}
\begin{equation}
\int_{0}^{\infty} dy \,
M^{(\nu, -)}_i(y; \xi_N) M^{(\nu, +)}_k(y; \xi_N)
=\delta_{ik}, \quad
0 \leq i, k \leq N-1.
\label{eqn:orth1}
\end{equation}

\begin{lem}
\label{thm:M+p}
Let $N \in \N, \xi_N \in \mM^+$ with $\xi_N(\R_+) =N$.
For $0 \leq s \leq t, x, y \in \R_+, 0 \leq i \leq N-1$,
\begin{eqnarray}
\label{eqn:Mplusp}
&& \int_{0}^{\infty} dx \,
p^{(\nu)}(t-s, y|x)
M^{(\nu, +)}_i \left( \frac{x}{s}; 
\frac{1}{s} \circ \xi_N \right)
= \left( \frac{s}{t} \right)^{i+1}
M^{(\nu, +)}_i \left( \frac{y}{t};
\frac{1}{t} \circ \xi_N \right) \\
\label{eqn:Mminusp}
&& \int_{0}^{\infty} dy \,
M^{(\nu, -)}_i \left( \frac{y}{t}; 
\frac{1}{t} \circ \xi_N \right)
p^{(\nu)}(t-s, y|x)
= \left( \frac{s}{t} \right)^{i}
M^{(\nu, -)}_i \left( \frac{x}{s};
\frac{1}{s} \circ \xi_N \right).
\end{eqnarray}
\end{lem}
\noindent{\it Proof.} \quad
By definition (\ref{eqn:M+-def}) and 
Lemma \ref{thm:QPintegral}
\begin{eqnarray}
&& M_i^{(\nu,+)}\left( \frac{x}{s}; \frac{1}{s} \circ \xi_N \right)
\nonumber\\
&=& \frac{1}{2 \pi \sqrt{-1}} 
\oint_{\Gamma((1/s) \circ \xi_{N, i+1})} dz \,
\frac{1}{2} \frac{1}{z^{\nu/2}} 
\left( \frac{x}{s} \right)^{\nu/2} 
I_{\nu} \left( \sqrt{z} \sqrt{\frac{x}{s}} \right)
e^{-z/2-x/(2s)}
\frac{1}{\prod_{a \in \xi_{N,i+1}} (z-a/s)}.
\nonumber
\end{eqnarray}
By setting $z=w/s$ and 
by the definition (\ref{eqn:pnu-}) of $p^{(\nu)}$, we find
\begin{eqnarray}
&& M_i^{(\nu, +)} \left( \frac{x}{s}; 
\frac{1}{s} \circ \xi_N \right)
\nonumber\\
&=& s^{i+1} \frac{1}{2 \pi \sqrt{-1}}
\oint_{\Gamma(\xi_{N, i+1})}dw \,
\frac{1}{2s} \left( \frac{x}{w} \right)^{\nu/2}
I_{\nu} \left( \frac{\sqrt{w x}}{s} \right)
e^{-(w+x)/(2s)}
\frac{1}{
\prod_{a \in \xi_{N,i+1}}(w-a)}
\nonumber\\
&=& s^{i+1} \frac{1}{2 \pi \sqrt{-1}}
\oint_{\Gamma(\xi_{N,i+1})} dw \,
p^{(\nu)}(s, x|w) 
\frac{1}{\prod_{a \in \xi_{N, i+1}}(w-a)}.
\label{eqn:M+rep}
\end{eqnarray}
Therefore the LHS of (\ref{eqn:Mplusp}) is 
$$
s^{i+1} \oint_{\Gamma(\xi_{N,i+1})} dw \,
\frac{1}{\prod_{a \in \xi_{N,i+1}} (w-a)}
\int_{0}^{\infty} dx \,
p^{(\nu)}(t-s, y|x) p^{(\nu)}(s, x|w).
$$
By the Chapman-Kolmogorov equation (\ref{eqn:CK1}),
it equals to
$$
s^{i+1} \frac{1}{2 \pi \sqrt{-1}}
\oint_{\Gamma(\xi_{N,i+1})} dw \,
\frac{1}{\prod_{a \in \xi_{N,i+1}}(w-a)}
p^{(\nu)}(t, y|w).
$$
Then, comparing with the expression (\ref{eqn:M+rep}),
we have (\ref{eqn:Mplusp}).

By definition (\ref{eqn:M+-def}) and Lemma \ref{thm:QPintegral}
$$
M^{(\nu, -)}_i \left( \frac{y}{t}; 
\frac{1}{t} \circ \xi_N \right)
= \int_{-\infty}^0 dw \,
\frac{1}{2} w^{\nu/2} \left( \frac{t}{y} \right)^{\nu/2}
e^{y/2t} e^{w/2} I_{\nu} \left( - \sqrt{w} \sqrt{\frac{y}{t}} \right)
\prod_{a \in (1/t) \circ \xi_{N,i}} (w-a).
$$
By setting $w=u/t$ and 
by the definition (\ref{eqn:pnu-}) of $p^{(\nu)}$, 
we have the expression
\begin{equation}
M^{(\nu, -)}_i \left( \frac{y}{t}; 
\frac{1}{t} \circ \xi_N \right)
=t^{-i} \int_{-\infty}^{0} du \,
p^{(\nu)}(-t, u|y) \prod_{a \in \xi_{N,i}}(u-a).
\label{eqn:M-rep}
\end{equation}
Therefore the LHS of (\ref{eqn:Mminusp}) is
$$
t^{-i} \int_{-\infty}^{0} du \,
\prod_{a \in \xi_{N,i}}(u-a)
\int_{0}^{\infty} dy \, 
p^{(\nu)}(-t, u|y) p^{(\nu)}(t-s, y|x).
$$
By the extension of the Chapman-Kolmogorov equation
(\ref{eqn:pnu-2}), it equals to
$$
t^{-i} \int_{-\infty}^{0} du \,
\prod_{a \in \xi_{N,i}}(u-a) p^{(\nu)}(-s, u|x)
= \left( \frac{s}{t} \right)^i
M^{(\nu, -)}_i \left( \frac{x}{s}; 
\frac{1}{s} \circ \xi_N \right),
$$
where the expression (\ref{eqn:M-rep}) was used.
This completes the proof of (\ref{eqn:Mminusp}).
\qed
\vskip 0.3cm

\begin{lem}
\label{thm:detM+}
Let $\y=(y_i)_{i=1}^{N} \in \W_N^{+}$.
For any 
$\xi_N (\cdot) = \sum_{i=1}^{N} \delta_{x_i}(\cdot)
\in \mM^+$ with a labeled configuration
$\x=(x_i)_{i=1}^{N}$ 
such that $x_1 \leq x_2 \leq \cdots \leq x_N$,
\begin{equation}
\frac{1}{h_N(\x)}
\det_{1 \leq i, j \leq N} \left[
\frac{1}{2} \left(\frac{y_j}{x_i} \right)^{\nu/2}
e^{-(x_i+y_j)/2} I_{\nu}(\sqrt{x_i y_j})
\right]
=\det_{1 \leq i, j \leq N}
\Big[ M^{(\nu, +)}_{i-1}(y_j; \xi_N) \Big].
\label{eqn:detM+}
\end{equation}
Here, when some of the $x_i$'s coincide, 
we interpret the LHS using 
l'H\^opital's rule.
\end{lem}
\noindent{\it Proof.} 
By the multilinearity of determinant
\begin{eqnarray}
&& \frac{1}{h_N(\x)}
\det_{1 \leq i, j \leq N} \left[
\frac{1}{2} \left(\frac{y_j}{x_i} \right)^{\nu/2}
e^{-(x_i+y_j)/2} I_{\nu}(\sqrt{x_i y_j})
\right]
\nonumber\\
&=& \det_{1 \leq i, j \leq N}
\left[ 
\frac{1}{2} \left(\frac{y_j}{x_i} \right)^{\nu/2}
e^{-(x_i+y_j)/2} I_{\nu}(\sqrt{x_i y_j})
\frac{1}{\prod_{k=1}^{i-1} (x_i-x_k)} \right]
\nonumber\\
&=& \det_{1 \leq i, j \leq N}
\left[ \sum_{\ell=1}^{i} 
\frac{1}{2} \left(\frac{y_j}{x_{\ell}} \right)^{\nu/2}
e^{-(x_{\ell}+y_j)/2} I_{\nu}(\sqrt{x_{\ell} y_j})
\frac{1}{\prod_{1 \leq k \leq i, k \not=\ell}
(x_{\ell}-x_k)} \right].
\nonumber
\end{eqnarray}
By definition (\ref{eqn:M+-def}) with (\ref{eqn:Qint})
of Lemma \ref{thm:QPintegral},
when $\xi_N \in \mM_0^+, \xi_N(\R_+)=N$,
\begin{eqnarray}
M^{(\nu, +)}_{i-1}(y; \xi_N) &=& 
\frac{1}{2 \pi \sqrt{-1}}
\oint_{\Gamma(\xi_{N,i})} d z \,
\frac{1}{2} \left(\frac{y}{z} \right)^{\nu/2}
e^{-(z+y)/2} I_{\nu}(\sqrt{z y})
\frac{1}{\prod_{x \in \xi_{N,i}}(z-x)} \nonumber\\
&=& 
\frac{1}{2 \pi \sqrt{-1}}
\oint_{\Gamma(\xi_{N,i})} d z \,
\frac{1}{2} \left(\frac{y}{z} \right)^{\nu/2}
e^{-(z+y)/2} I_{\nu}(\sqrt{z y})
\frac{1}{\prod_{k=1}^{i} (z-x_{k})} 
\nonumber\\
&=& \sum_{\ell=1}^{i} 
\frac{1}{2} \left(\frac{y}{x_{\ell}} \right)^{\nu/2}
e^{-(x_{\ell}+y)/2} I_{\nu}(\sqrt{x_{\ell} y})
\frac{1}{\prod_{1 \leq k \leq i, k \not= \ell}
(x_{\ell}-x_k)},
\label{eqn:M+i-1A}
\end{eqnarray}
$1 \leq i \leq N, y \in \R_+$.
Then (\ref{eqn:detM+}) is proved for $\xi_N \in \mM_0^+$. 
When some of the $x_i$'s coincide,
the LHS of (\ref{eqn:detM+}) is interpreted
using l'H\^opital's rule and
in the RHS of (\ref{eqn:detM+})
$M^{(\nu,+)}_{i-1}(y_j; \xi_N)$ should be given by
the first expression of (\ref{eqn:M+i-1A}).
Then (\ref{eqn:detM+}) is valid for any
$\xi_N \in \mM^+$ with $\xi_N(\R_+)=N \in \N$. \qed
\vskip 0.3cm

\SSC{Proofs of Theorems}
\subsection{Proof of Theorem \ref{thm:Finite}}
In this subsection we give a proof of 
Theorem \ref{thm:Finite}.
Assume that $\xi_N \in \mM^+$ with $\xi_N(\R_+) = N \in \N$.
Define 
\begin{eqnarray}
\label{eqn:phin+}
\phi_{i}^{(\nu, +)}(t, x; \xi_N)
&\equiv& 
t^{-(i+1)}
M_{i}^{(\nu, +)} \left( \frac{x}{t}; 
\frac{1}{t} \circ \xi_N \right), \\
\label{eqn:phin-}
\phi_{i}^{(\nu, -)}(t,x; \xi_N)
&\equiv& t^{i}
M_{i}^{(\nu, -)}\left( \frac{x}{t}; 
\frac{1}{t} \circ \xi_N \right),
\end{eqnarray}
$0 \leq i \leq N-1, t>0, x \in \R$. 
From Lemma \ref{thm:M+p}, 
the following relations
are derived.

\begin{lem}
\label{thm:phi+-p}
For $\xi_N \in \mM_0$ with $\xi_N(\R)=N \in \N$,
$0 \leq t_1 \leq t_2$,
\begin{eqnarray}
\label{eqn:key1}
&& \int_{0}^{\infty} dx_2 \,
\phi_{i}^{(\nu,-)}(t_2, x_2; \xi_N) p^{(\nu)}(t_2-t_1, x_2|x_1)
=\phi_{i}^{(\nu,-)}(t_1, x_1; \xi_N), \quad
0 \leq i \leq N-1, \quad \\
\label{eqn:key2}
&& \int_{0}^{\infty} dx_1 \,
p^{(\nu)}(t_2-t_1, x_2 |x_1)
\phi_{i}^{(\nu,+)}(t_1, x_1; \xi_N)
=\phi_{i}^{(\nu,+)}(t_2, x_2; \xi_N), \quad
0 \leq i \leq N-1, \quad \\
&& \int_{0}^{\infty} dx_1 
\int_{0}^{\infty} dx_2 \,
\phi_{i}^{(\nu,-)}(t_2, x_2; \xi_N)
p^{(\nu)}(t_2-t_1, x_2|x_1)
\phi_{j}^{(\nu,+)}(t_1, x_1; \xi_N)
=\delta_{ij}, 
\nonumber\\
\label{eqn:key3}
&& \hskip 10cm
0 \leq i, j \leq N-1. 
\end{eqnarray}
\end{lem}
\vskip 0.3cm
\noindent{\bf Remark 3.} 
The equation (\ref{eqn:key3}) is obtained
from the combination of (\ref{eqn:key1}) and (\ref{eqn:key2}).
It should be emphasize the fact that
the biorthonormality (\ref{eqn:orth1})
is obtained by just taking the limit
$t_2-t_1 \to 0$ in (\ref{eqn:key3}).
For the purpose in this paper,
the connection to the theory of
multiple orthogonal polynomials
\cite{CA03a,CA03b,BK05,KMW09,DF08} is not so essential.
We can think that the integrals (\ref{eqn:M+rep})
and (\ref{eqn:M-rep}) define the functions
$\phi_i^{(\nu,+)}(t,x;\xi)$
and $\phi_i^{(\nu,-)}(t,x;\xi)$ 
through (\ref{eqn:phin+}) and
(\ref{eqn:phin-}), respectively.
As shown in the proof of Lemma \ref{thm:M+p},
(\ref{eqn:key1}) and (\ref{eqn:key2}) are
readily obtained by applying the 
Chapman-Kolmogorov equation (\ref{eqn:CK1})
and its extension (\ref{eqn:pnu-2}).
\vskip 0.3cm

By definitions (\ref{eqn:type2a}), (\ref{eqn:M+-def}), and
(\ref{eqn:phin-}),
we can see that $\phi_i^{(\nu,-)}(t,x; \xi_N)$
is a monic polynomial of $x$ of degree $i$, which is 
independent of $t$.
By this fact and Lemma \ref{thm:detM+}, 
the multitime probability density 
(\ref{eqn:pxi})
is expressed as
\begin{eqnarray}
&& p^{\xi_N}_{\nu} \Big(t_1, \xi_N^{(1)}; \dots ;
t_M, \xi_N^{(M)} \Big) 
= \det_{1 \leq i, j \leq N}
\Big[ \phi^{(\nu,-)}_{i-1}(t_M, x^{(M)}_j; \xi_N) \Big]
\nonumber\\
&& \qquad \times
\prod_{m=1}^{M-1} f_N^{(\nu)}(t_{m+1}-t_m; \x^{(m+1)}|\x^{(m)})
\det_{1 \leq k, \ell \leq N}
\Big[ \phi^{(\nu,+)}_{k-1}(t_1, x^{(1)}_{\ell}; \xi_N) \Big]
\label{eqn:pnuxiN2}
\end{eqnarray}
for $\xi_N \in \mM^+$ with $\xi_N(\R_+)=N \in \N$.
By the argument given in Section 4 in \cite{KT07b},
the expression (\ref{eqn:pnuxiN2}) with
Lemma \ref{thm:phi+-p} leads to
the Fredholm determinantal expression
for the generating function
of multitime correlation functions (\ref{eqn:G1}),
$$
\cG^{\xi_N}_{\nu}[\vchi]=
\Det_{
\substack{1 \leq m, n \leq M \\ (x,y) \in (0, \infty)^2}
}
\Bigg[ \delta_{m n} \delta (x-y)
+\widetilde{S}^{m,n}(x,y; \xi_N)
\chi_n (y) \Bigg],
$$
where
$$
\widetilde{S}^{m,n}(x,y; \xi_N)
=S^{m,n}(x,y; \xi_N)
-{\bf 1}(m > n) p^{(\nu)}(t_m-t_n, x|y)
$$
with
\begin{eqnarray}
&& S^{m,n}(x,y; \xi_N)
= \sum_{i=0}^{N-1} 
\phi_{i}^{(\nu, +)}(t_m, x; \xi_N)
\phi_{i}^{(\nu, -)}(t_n, y; \xi_N)
\nonumber\\
&& \qquad =\frac{1}{t_m}
\sum_{i=0}^{N-1}
\left( \frac{t_n}{t_m} \right)^{i/2}
M_i^{(\nu,+)}\left( \frac{x}{t_m}; 
\frac{1}{t_m} \circ \xi_N \right)
M_i^{(\nu,-)} \left( \frac{y}{t_n};
\frac{1}{t_n} \circ \xi_N \right).
\label{eqn:Smn0}
\end{eqnarray}
Here the Fredholm determinant is defined by
the following expansion
\begin{eqnarray}
&& \Det_{
\substack{1 \leq m, n \leq M \\ (x,y) \in (0, \infty)^2}
}
 \Bigg[
\delta_{m n} \delta(x-y)
+\widetilde{S}^{m,n}(x,y; \xi_N) \chi_{n}(y) \Bigg]
\nonumber\\
&=& \sum_{N_{1}=0}^{N} \cdots
\sum_{N_{M}=0}^{N}
\prod_{m=1}^{M}\frac{1}{N_m !}
\int_{\R_+^{N_{1}}} \prod_{i=1}^{N_1} d x_{i}^{(1)}
 \cdots
\int_{\R_+^{N_{M}}} 
\prod_{i=1}^{N_{M}} d x_{i}^{(M)} \nonumber\\
&& \quad \quad
\times \prod_{m=1}^{M} \prod_{i=1}^{N_{m}} 
\chi_{m} \Big(x_{i}^{(m)} \Big) 
\det_{
\substack{1 \leq i \leq N_{m}, 1 \leq j \leq N_{n} \\ 1 \leq m, n \leq M}
}
 \Bigg[
\widetilde{S}^{m,n}(x^{(m)}_i, x^{(n)}_j; \xi_N)
\Bigg].
\label{eqn:Fredholm2}
\end{eqnarray}

The following invariance of 
finite dimensional distributions of determinantal 
processes will be used.

\begin{lem}\label{thm:gauge}
Let $(\Xi(t), \P^{\xi})$ and $(\widetilde{\Xi}(t), \widetilde{\P}^{\xi})$
be the $\mM^+$-valued processes, which are determinantal with
correlation kernels $\mbK$ and $\widetilde{\mbK}$,
respectively.
If there is a function $G(s,x)$, which is continuous 
with respect to $x \in (0, \infty)$ for any fixed $s \in [0, \infty)$,
such that
\begin{equation}
\mbK(s,x;t,y)=\frac{G(s,x)}{G(t,y)}
\widetilde{\mbK}(s,x;t,y),
\quad (s, t) \in [0, \infty)^2, \quad (x,y) \in (0, \infty)^2,
\label{eqn:gauge1}
\end{equation}
then
\begin{equation}
(\Xi(t), \P^{\xi})=(\widetilde{\Xi}(t), \widetilde{\P}^{\xi})
\label{eqn:gauge2}
\end{equation}
in the sense of finite dimensional distributions.
\end{lem}
\vskip 0.3cm

The relation (\ref{eqn:gauge1})
is called the {\it gauge transformation} and
(\ref{eqn:gauge2}) is said to be the {\it gauge invariance}
of the determinantal processes \cite{KT09}.

\vskip 0.3cm
\noindent{\it Proof of Theorem \ref{thm:Finite}.} 
Inserting the integral representations
for $M_i^{(\nu, \pm)}$, 
given by (\ref{eqn:M+rep}) and (\ref{eqn:M-rep}), 
into (\ref{eqn:Smn0}), 
the kernel $S^{m,n}$ is written as 
\begin{eqnarray}
S^{m,n}(x,y; \xi_N) 
&=& \frac{1}{2\pi \sqrt{-1}} 
\int_{-\infty}^{0} du \,
\oint_{\Gamma(t_m^{-1} \circ \xi_N)} dz \,
p^{(\nu)}(t_m, x|z) p^{(\nu)}(-t_n, u|y)
\nonumber\\
&& \qquad \qquad \times
\sum_{k=0}^{N-1}
\frac{\prod_{j=1}^{i}(u-a_{j})}
{\prod_{j=1}^{i+1} (z-a_{j})},
\nonumber
\end{eqnarray}
where $\xi_N(\cdot)=\sum_{i=1}^{N} \delta_{a_i}(\cdot)$.
For $z_1, z_2 \in \C$
with $z_1 \notin \{x_1, \dots, x_N\}$, the following identity
holds,
$$
\sum_{k=0}^{N-1} 
\frac{\prod_{j=1}^{k} (z_2-x_{j})}
{\prod_{j=1}^{k+1} (z_1-x_{j})}
= \left( \prod_{j=1}^{N}
\frac{z_2-x_{j}}{z_1-x_{j}}
-1 \right)
\frac{1}{z_2-z_1}.
$$
By this identity, we have
\begin{eqnarray}
S^{m,n}(x,y;\xi_N) 
&=& 
 \frac{1}{2\pi \sqrt{-1}} 
\int_{-\infty}^{0} du \,
\oint_{\Gamma_{u}(t_m^{-1} \circ \xi_N)} dz \,
p^{(\nu)}(t_m, x|z) p^{(\nu)}(-t_n, u|y)
\nonumber\\
&& \qquad \times
\frac{1}{u-z}
\left( \prod_{j=1}^{N}
\frac{u-a_{j}}{z-a_{j}}
-1 \right).
\nonumber
\end{eqnarray}
Note that
$$
 \frac{1}{2\pi \sqrt{-1}} 
\int_{-\infty}^{-\varepsilon} du \,
\oint_{\Gamma_{u}(t_m^{-1} \circ \xi_N)} dz \,
p^{(\nu)}(t_m, x|z) p^{(\nu)}(-t_n,u|y)
\frac{1}{u-z}=0, \quad
\forall \varepsilon >0,
$$
since by definition (\ref{eqn:pnu-})
$$
p^{(\nu)}(t_m, x|z)
=\frac{x^{\nu}}{(2t_{m})^{\nu+1}}
e^{-(x+z)/2t_m}
\sum_{n=0}^{\infty} \frac{1}{\Gamma(n+1) \Gamma(n+1+\nu)}
\frac{x^n z^n}{(2t_m)^{2n}}
$$
is an entire function with respect to $z$.
Then (\ref{eqn:KN1a}) is obtained.
When $\xi_N \in \mM^+_0$, the Cauchy integrals are performed
in (\ref{eqn:KN1a}) and (\ref{eqn:KN1bz}) is obtained.
If we use (\ref{eqn:pnupJ})
and (\ref{eqn:entire3}), 
we have the equality
\begin{equation}
\mbK_{\nu}^{\xi_N}(s, x; t, y)
= \left(\frac{x}{y} \right)^{\nu/2}
\mbK^{\xi_N}_{J_{\nu}}(s, x; t, y),
\quad (s,t) \in [0, \infty)^2, (x,y) \in (0, \infty)^2.
\label{eqn:gauge3}
\end{equation}

Then, by Lemma \ref{thm:gauge}, 
the correlation kernel $\mbK_{\nu}^{\xi_N}$
can be replaced by $\mbK^{\xi_N}_{J_{\nu}}$
without changing any finite dimensional distributions.
This completes the proof.
\qed
\vskip 0.3cm

\subsection{Proof of Theorem \ref{thm:Infinite}}

In this subsection we give a proof of Theorem \ref{thm:Infinite}.
For $\zeta \in \mM, L > 0$, we define
$$
M(\zeta, L)=\int_{[-L, L] \setminus \{0\}}
\frac{\zeta(dx)}{x}
$$
and put
$$
M(\zeta)=\lim_{L \to \infty} M(\xi, L),
$$
if the limit finitely exists.
In an earlier paper we proved the following lemma
(Lemma 4.4 in \cite{KT10}).

\begin{lem}
\label{thm:Dyson}
Assume that $\zeta \in \mM$ satisfies the following conditions: \\
\noindent {\rm {\bf (C.1)}}
there exists $\widehat{C}_0 > 0$ such that
$|M(\zeta)|  < \widehat{C}_0$, \\
\noindent {\rm {\bf (C.2)} (i) }
there exist $\widehat{\alpha} \in (1,2)$ and 
$\widehat{C}_1>0$ such that
$M_{\widehat{\alpha}}(\zeta) \le \widehat{C}_1$, 
\noindent {\rm (ii)}
there exist $\widehat{\beta} >0$ and $\widehat{C}_2 >0$ such that
$$
M_1(\tau_{-a^2} \zeta^{\langle 2 \rangle}) \leq \widehat{C}_2
(|a| \vee 1)^{-\widehat{\beta}},
\quad \forall a \in \supp \zeta.
$$
Then there exists 
$\widehat{C}_3=\widehat{C}_3(\widehat{\alpha}, 
\widehat{\beta}, \widehat{C}_0,
\widehat{C}_1,\widehat{C}_2) >0$
and $\theta \in (\widehat{\alpha} 
\vee (2-\widehat{\beta}), 2)$ such that
$$
|\Phi_0(\zeta, a, \sqrt{-1} y)|
\leq \Big[ \widehat{C}_3
\Big\{ (|y| \vee 1)^{\widehat{\theta}}
+(|a| \vee 1)^{\widehat{\theta}} \Big\} \Big],
\quad \forall y \in \R,
\quad \forall a \in \supp \zeta.
$$
\end{lem}

Here we use the following version.

\begin{lem}
\label{thm:Pi_bound}
For any $\xi \in \mX_0^+$, there exist
$C_3=C_3(\alpha, \beta, C_1, C_2) >0$
and $\theta \in (\alpha \vee (1-\beta), 1)$
such that
$$
|\Phi_0(\xi, a, y)|
\leq \exp \Big[ C_3
\Big\{ (|a| \vee 1)^{\theta} 
+(|y| \vee 1)^{\theta} \Big\} \Big],
\quad \forall y <0,
\quad \forall a \in \supp \xi.
$$
\end{lem}
\noindent{\it Proof.} \,
First we consider the case that $a=0 \in \xi$.
Let $\zeta=\xi^{\langle 1/2 \rangle}$. Then
the equalities
\begin{eqnarray}
\Phi_0(\xi, 0, z)
&=& \prod_{x \in \xi \cap \{0\}^{\rm c}}
\left( 1- \frac{z}{x} \right)
= \prod_{x \in \xi \cap \{0\}^{\rm c}}
\left( 1- \frac{\sqrt{z}}{\sqrt{x}} \right)
\left( 1 + \frac{\sqrt{z}}{\sqrt{x}} \right)
\nonumber\\
&=& \prod_{x' \in \zeta \cap \{0\}^{\rm c}}
\left( 1- \frac{\sqrt{z}}{x'} \right)
= \Phi_{0}(\zeta, 0, \sqrt{z})
\label{eqn:eqPhi}
\end{eqnarray}
hold, if the products finitely exist.
Since $M(\zeta)=0$ by the definition
$\zeta=\xi^{\langle 1/2 \rangle}$, 
$\zeta$ satisfies the condition {\bf (C.1)}
in Lemma \ref{thm:Dyson}.
By assumption of the present lemma, 
that is, $\xi$ satisfies the condition {\bf (C.A)},
the condition {\bf (C.2)}
is satisfied by $\zeta$ with
$\widehat{\alpha}=2 \alpha$ and $\widehat{\beta}=2 \beta$,
since
$$
M_{\alpha}(\xi)=(M_{2 \alpha}(\zeta))^2
\quad \mbox{and} \quad
M_1(\tau_{-a} \xi)=M_1(\tau_{-(\sqrt{a})^2} 
\zeta^{\langle 2 \rangle}),
\quad a \in \supp \xi.
$$
Then, by Lemma \ref{thm:Dyson}, 
there exist
$\widehat{\theta}
\in (\widehat{\alpha} \vee (2-\widehat{\beta}), 2)
=(2 \alpha \vee 2(1-\beta), 2)$ and
$\widehat{C}_3 < \infty$ such that
$
\Phi_0(\zeta, 0, \sqrt{-z}) \leq 
\exp \Big[ \widehat{C}_3 (\sqrt{z} \vee 1)^{\widehat{\theta}} \Big],
z>0$.
Since (\ref{eqn:eqPhi}) is established,
it implies
\begin{equation}
\Phi_0(\xi, 0, -z) \leq 
\exp \Big[ \widehat{C}_3 (|z| \vee 1)^{\theta} \Big]
\label{eqn:estimate1}
\end{equation}
with $\theta=\widehat{\theta}/2$.

Next we consider the case that $a \in \supp \xi$
and $a \not=0$.
For $y < 0$, the equality
$$
\Phi_0(\xi, a, y)
= \Phi_{0}(\xi, 0, y) \Phi_{0}(\xi \cap \{0\}^{\rm c}, a, 0)
\left( \frac{y}{a} \right)^{\xi(\{0\})}
\frac{a}{a-y}
$$
is valid. Since $y < 0$ is assumed
$$
\left| \left( \frac{y}{a} \right)^{\xi(\{0\})}
\frac{a}{a-y} \right|
=\left| \frac{a}{a-y} \right|
\vee \left| \frac{y}{a-y} \right| 
\leq 1.
$$
On the other hand,
\begin{eqnarray}
 \Phi_0(\xi \cap \{0\}^{\rm c}, a, 0)
&=& \prod_{a \in \xi \cap \{0\}^{\rm c}}
\left( 1+ \frac{a}{x-a} \right)
\nonumber\\
&\leq& \exp \left\{ \int_{\R}
(\tau_{-a} \xi) (dx) \frac{|a|}{|x|} \right\}
= \exp \Big\{ |a| M_1(\tau_{-a} \xi) \Big\}.
\nonumber
\end{eqnarray}
By the condition {\bf (C.A)} (ii) of $\mX_0^+$,
it is bounded from the above by
$
\exp \{ C_2 (|a| \vee 1)^{-\beta} |a| \}
= \exp \{ C_2 (|a| \vee 1)^{1-\beta} \}.
$
Combining this with (\ref{eqn:estimate1}),
the proof is completed. \qed
\vskip 0.3cm

\noindent{\it Proof of Theorem \ref{thm:Infinite}.} 
Note that $\xi\cap [0, L]$, $L>0$ and $\xi$ 
satisfy {\bf (C.A)} 
with the same constants $C_1, C_2$ and indices $\alpha, \beta$.
By Lemma \ref{thm:Pi_bound}
we see that there exists $C_3 >0$ such that
$$
|\Phi_0 (\xi \cap [0, L], x', y')| \leq
\exp \Big[ C_3 \Big\{(|x'|\vee1)^{\theta} 
+ (|y'|^{\theta} \vee 1) \Big\} \Big],
$$
$\forall L>0, \; 
\forall x' \in \supp \xi \cap [0, L], 
\forall y' < 0$
with $\theta \in (\alpha \vee (1-\beta), 1)$.
Therefore, since for any $x' \in \supp \xi$, $y' < 0$
$$
\Phi_0(\xi\cap [0, L], x', y') \to \Phi_0(\xi, x',  y'), 
\quad L\to\infty,
$$
we can apply Lebesgue's convergence theorem 
to (\ref{eqn:KN1b}) and obtain
$$
\lim_{L\to\infty}\mbK_{\nu}^{\xi\cap [0, L]}
\left(s, x; t, y \right)
=\mbK_{\nu}^{\xi}\left(s, x; t, y\right).
$$
Since for any $(s,t) \in (0, \infty)^{2}$ and 
any compact interval $I \subset (0, \infty)$
$$
\sup_{x, y \in I} \Big|
\mbK_{\nu}^{\xi \cap [0, L]}(s, x; t, y) \Big| < \infty,
$$
we can obtain the convergence of generating functions
for multitime correlation functions,
${\cal G}_{\nu}^{\xi \cap [0, L]}[\vchi] \to
{\cal G}_{\nu}^{\xi}[\vchi]$, as $L \to \infty$.
It implies 
$\P_{\nu}^{\xi \cap [0, L]} \to \P_{\nu}^{\xi}$ as $L \to \infty$
in the sense of finite dimensional distributions.
By virtue of the gauge invariance of
determinantal processes (Lemma \ref{thm:gauge}), 
the proof is completed.
\qed

\subsection{Proof of Theorem \ref{thm:RelBessel}}

In this subsection we give a proof of 
Theorem \ref{thm:RelBessel}.
First we prove a lemma.

\begin{lem}
\label{thm:Phinu}
For $z \not= j_{\nu, i}, i \in \N$
\begin{equation}
\Phi^{(\nu)}(\xi_{J_{\nu}}^{\langle 2 \rangle}, 
(j_{\nu, i})^2, z^2)
=\frac{2 j_{\nu, i}}{(j_{\nu, i})^2-z^2}
\frac{J_{\nu}(z)}
{J_{\nu+1}(j_{\nu, i})},
\label{eqn:Equal2}
\end{equation}
\end{lem}
\noindent{\it Proof.} \,
By the expression (\ref{eqn:entire1}) of $J_{\nu}(z)$
and the definition (\ref{eqn:entire2}) of
$\Pi^{(\nu)}$, we have
\begin{equation}
J_{\nu}(z)= \frac{2^{-\nu}}{\Gamma(\nu+1)} 
\Pi^{(\nu)}(\xi_{J_{\nu}}^{\langle 2 \rangle}, z^2).
\label{eqn:J_Phi1}
\end{equation}
Put 
$
{\Pi^{(\nu)}}'(\xi^{\langle 2 \rangle}, z^2)
\equiv (d/dz) \Pi^{(\nu)}(\xi^{\langle 2 \rangle}, z^2).
$
Then
\begin{equation}
J'_{\nu}(j_{\nu, i})= \frac{2^{-\nu}}{\Gamma(\nu+1)}
{\Pi^{(\nu)}}'
(\xi_{J_{\nu}}^{\langle 2 \rangle}, j_{\nu, i}^2).
\label{eqn:J_Phi2}
\end{equation}
We see that, if $x' \in \xi$,
$$
\frac{1}{1-z^2/(x')^2} \Pi^{(\nu)}(\xi^{\langle 2 \rangle}, z^2)
= z^{\nu}
\prod_{x \in \xi \cap \{x'\}^{\rm c}}
\left( 1- \frac{z^2}{x^2} \right).
$$
Since
$$
{\Pi^{(\nu)}}'(\xi^{\langle 2 \rangle}, z^2)
=\nu \frac{1}{z} \Pi^{(\nu)}(\xi^{\langle 2 \rangle}, z^2)
+z^{\nu}
\int_{\R_+} \xi(dx)
\left(- \frac{2z}{x^2} \right)
\prod_{y \in \xi \cap \{x\}^{\rm c}}
\left(1- \frac{z^2}{y^2} \right),
$$
if $x' \in \xi$, 
$$
{\Pi^{(\nu)}}'(\xi^{\langle 2 \rangle}, (x')^2)
=-\frac{2}{x'}
(x')^{\nu}
\prod_{x \in \xi \cap \{x'\}^{\rm c}}
\left(1-\frac{(x')^2}{x^2} \right).
$$
Therefore, for $ x' \in \xi$,
\begin{equation}
\frac{1}{1-z^2/(x')^2}
\frac{\Pi^{(\nu)}(\xi^{\langle 2 \rangle}, z^2)}
{(\Pi^{(\nu)})'(\xi^{\langle 2 \rangle}, (x')^2)}
=-\frac{x'}{2} \left( \frac{z}{x'} \right)^{\nu}
\prod_{x \in \xi \cap \{x'\}^{\rm c}}
\frac{1-z^2/x^2}{1-(x')^2/x^2}.
\label{eqn:eqE1}
\end{equation}
Since
$$
1-\frac{z^2-a^2}{x^2-a^2}
=\frac{x^2-z^2}{x^2-a^2}
=\frac{1-z^2/x^2}{1-a^2/x^2}
$$
for $x \not=a$, 
combining (\ref{eqn:eqE1}) with
the definition (\ref{eqn:entire3}) gives
\begin{eqnarray}
\Phi^{(\nu)}(\xi^{\langle 2 \rangle}, (x')^2, z^2)
&=& 
\left( \frac{z}{x'} \right)^{\nu} 
\prod_{x \in \xi \cap \{x'\}^{\rm c}}
\frac{ 1-z^2/x^2}{1-(x')^2/x^2}
\nonumber\\
&=&
\frac{2x'}{z^2-(x')^2}
\frac{\Pi^{(\nu)}(\xi^{\langle 2 \rangle}, z^2)}
{{\Pi^{(\nu)}}'(\xi^{\langle 2 \rangle}, (x')^2)},
\quad x' \in \xi. 
\label{eqn:Equal1}
\end{eqnarray}
Setting $\xi=\xi_{J_{\nu}}$ and $x'=j_{\nu, i}$
in (\ref{eqn:Equal1}), we have
$$
\Phi^{(\nu)}(\xi_{J_{\nu}}^{\langle 2 \rangle}, 
(j_{\nu, i})^2, z^2)
=
\frac{2 j_{\nu, i}}{z^2-(j_{\nu, i})^2}
\frac{J_{\nu}(z)}
{J'_{\nu}(j_{\nu, i})},
$$
by (\ref{eqn:J_Phi1}) and (\ref{eqn:J_Phi2}).
If we use (\ref{eqn:J'2}), (\ref{eqn:Equal2}) is obtained.
\qed
\vskip 0.3cm

\noindent{\it Proof of Theorem \ref{thm:RelBessel}.} \,
We set $\xi=\xi_{J_{\nu}}^{\langle 2 \rangle}$
in $\mbK_{J_{\nu}}^{\xi}$ and obtain the kernel
\begin{eqnarray}
\mbK_{J_{\nu}}(s,x; t, y)
&=& \sum_{i=1}^{\infty} \int_{-\infty}^{0} dz \,
p_{J_{\nu}}(s, x| j_{\nu, i}^2)
\Phi^{(\nu)}(\xi_{J_{\nu}}^{\langle 2 \rangle}, (j_{\nu, i})^2, z)
p_{J_{\nu}}(-t, z|y)
\nonumber\\
\label{eqn:AA1}
&& \quad - {\bf 1}(s>t)
p_{J_{\nu}}(s-t, x|y).
\end{eqnarray}
By Lemmas \ref{thm:j1} and \ref{thm:Phinu},
\begin{eqnarray}
\Phi^{(\nu)}(\xi_{J_{\nu}}^{\langle 2 \rangle},
(j_{\nu, i})^2, z)
&=& \frac{2 j_{\nu, i}}{(j_{\nu, i})^2-z}
\frac{J_{\nu}(\sqrt{z})}{J_{\nu+1}(j_{\nu, i})}
\nonumber\\
&=& \frac{1}{(J_{\nu+1}(j_{\nu, i}))^2}
\int_{0}^{1} du \,
J_{\nu}(\sqrt{uz}) J_{\nu}(\sqrt{u} j_{\nu, i}).
\nonumber
\end{eqnarray}
By the integral formula (\ref{eqn:pJnu-}),
the first term of the RHS of (\ref{eqn:AA1}) equals
$$
\sum_{i=1}^{\infty} p_{J_{\nu}}(s,x|(j_{\nu, i})^2)
\frac{1}{(J_{\nu+1}(j_{\nu, i}))^2}
\int_{0}^{1} du \, e^{ut/2}
J_{\nu}(\sqrt{u y}) J_{\nu}(\sqrt{u} j_{\nu, i}).
$$
Since $s > 0$, we can use the expression (\ref{eqn:pJnu1})
for $p_{J_{\nu}}(s,x|j_{\nu, i}^2)$
and the above is written as
$$
\int_{0}^{1} du 
\int_{0}^{\infty} d w \,
e^{ut/2-2ws} J_{\nu}(\sqrt{uy}) J_{\nu}(2 \sqrt{wx})
\sum_{i=1}^{\infty} 
\frac{J_{\nu}(2 \sqrt{w} j_{\nu, i})
J_{\nu}(\sqrt{u} j_{\nu, i})}
{(J_{\nu+1}(j_{\nu, i}))^2}.
$$
From the completeness (\ref{eqn:Fourier-Bessel})
of the orthonormal basis 
$\{J_{\nu}(j_{\nu,i} \sqrt{x})/J_{\nu+1}(j_{\nu,i}), i \in \N\}$
in $L^{2}(0,1)$,
the above gives
$$
\mbK_{J_{\nu}}(s,x;t,y)
=\bK_{J_{\nu}}(s,x;t,y)+R(s,x;t,y)
$$
with the extended Bessel kernel (\ref{eqn:KBessel1})
and
\begin{eqnarray}
R(s,x;t,y)
&=& \int_{0}^{1} du \int_{1}^{\infty} dw \,
e^{ut/2-2ws} J_{\nu}(\sqrt{uy})
J_{\nu}(2 \sqrt{w x})
\nonumber\\
&& \times
\sum_{i=1}^{\infty} 
\frac{J_{\nu}(2 \sqrt{w} j_{\nu, i})
J_{\nu}(\sqrt{u} j_{\nu, i})}
{(J_{\nu+1}(j_{\nu, i}))^2}.
\nonumber
\end{eqnarray}
Since for any fixed $s,t > 0$
$$
|R(s+\theta, x: t+\theta, y)| \to 0
\quad \mbox{as} \quad
\theta \to \infty
$$
uniformly on any compact subset of
$(x,y) \in \R_{+}^2$,
$$
\mbK_{J_{\nu}}(s+\theta, x; t+\theta)
\to
\bK_{J_{\nu}}(s,x; t, y)
\quad \mbox{as} \quad \theta \to \infty
$$
holds in the same sense.
Hence we obtain (\ref{eqn:relaxBessel}).
This completes the proof.
\qed
\vskip 0.3cm

\begin{small}
\noindent{\bf Acknowledgements} \quad
M.K. is supported in part by
the Grant-in-Aid for Scientific Research (C)
(No.21540397) of Japan Society for
the Promotion of Science.
H.T. is supported in part by
the Grant-in-Aid for Scientific Research 
(KIBAN-C, No.19540114) of Japan Society for
the Promotion of Science.


\end{small}
\end{document}